\documentclass[11pt]{amsart}
\usepackage[T1]{fontenc}
\usepackage{amsmath,amssymb,amsthm,mathrsfs}
\usepackage[a4paper,margin=25mm]{geometry}

\usepackage{enumitem}
\setlist[enumerate]{font=\upshape}
\usepackage{tikz-cd}
\usepackage[colorlinks=true,linkcolor=blue,citecolor=blue,urlcolor=blue,pagebackref]{hyperref}
\usepackage[capitalize,nameinlink]{cleveref}

\numberwithin{equation}{section}
\newtheorem{theorem}{Theorem}[section]
\newtheorem{proposition}[theorem]{Proposition}
\newtheorem{lemma}[theorem]{Lemma}
\newtheorem{corollary}[theorem]{Corollary}
\newtheorem{claim}[theorem]{Claim}
\theoremstyle{definition}
\newtheorem{definition}[theorem]{Definition}
\newtheorem{assumption}[theorem]{Assumption}
\newtheorem{conjecture}[theorem]{Conjecture}

\theoremstyle{remark}
\newtheorem{remark}[theorem]{Remark}

\newtheorem*{theorem*}{Theorem}

\AddToHook{env/proposition/begin}{\crefalias{theorem}{proposition}}
\AddToHook{env/lemma/begin}{\crefalias{theorem}{lemma}}
\AddToHook{env/corollary/begin}{\crefalias{theorem}{corollary}}
\AddToHook{env/claim/begin}{\crefalias{theorem}{claim}}
\AddToHook{env/definition/begin}{\crefalias{theorem}{definition}}
\AddToHook{env/assumption/begin}{\crefalias{theorem}{assumption}}
\AddToHook{env/conjecture/begin}{\crefalias{theorem}{conjecture}}
\AddToHook{env/notation/begin}{\crefalias{theorem}{notation}}
\AddToHook{env/remark/begin}{\crefalias{theorem}{remark}}
\AddToHook{env/example/begin}{\crefalias{theorem}{example}}

\crefname{conjecture}{Conjecture}{Conjectures}

\newcommand{\cX}{\mathcal{X}}

\newcommand{\cY}{\mathcal{Y}}

\newcommand{\cO}{\mathcal{O}}
\newcommand{\Q}{\mathbb Q}
\newcommand{\R}{\mathbb R}
\newcommand{\C}{\mathbb C}
\newcommand{\BC}{\mathrm{BC}}
\newcommand{\bbeta}{\boldsymbol{\beta}}

\DeclareMathOperator{\Supp}{Supp}

\DeclareMathOperator{\mult}{mult}
\DeclareMathOperator{\Exc}{Exc}

\title{Invariance of plurigenera under the good Zariski decomposition assumption}
\author{Yi Li}
 	\address{Department of Mathematics, Wuhan University, Hubei, Wuhan 430074, China}
 	\email{yilimath@whu.edu.cn}
 \author{Sheng Rao}
 	\address{Department of Mathematics, Wuhan University, Hubei, Wuhan 430074, China}
 	\email{likeanyone@whu.edu.cn}
\date{\today}

\begin{document}
	
	\begin{abstract}
		For a smooth family of compact K\"ahler manifolds over a disk, we prove that all plurigenera of the fibers are constant, assuming the good divisorial Zariski decomposition. This settles Siu's conjecture on deformation invariance of plurigenera for compact K\"ahler manifolds under that assumption. The main new ingredient is a termination result: on a compact K\"ahler manifold with pseudo-effective canonical class, a minimal model program terminates after finitely many divisorial contractions and flips, with a good minimal model. Following an idea of Koll\'ar and using the extension theorem of Hacon--Xie and Hacon--Li--Rao, we extend this  minimal model program to a family. Hence, by  Levine--Koll\'ar--Mori's deformation invariance of plurigenera for varieties in the Fujiki's class $\mathcal{C}$ with terminal singularities and semiample canonical diviors, the plurigenera are locally constant in this family and hence the original smooth fiberwise K\"ahler family.
	\end{abstract}
	
	\maketitle
	\setcounter{tocdepth}{1}
	\tableofcontents

	\section{Introduction}\label{sec:intro}

	The plurigenera are among the most fundamental birational invariants of a compact complex manifold. Their behaviour is at the heart of the classification theory of algebraic varieties and of the minimal model program. It is natural to ask how they behave when $X$ varies in a family. For smooth projective families this is answered by Siu's celebrated theorem on the invariance of plurigenera \cite{Si98,Siu02}, one of the deepest applications of $L^2$ methods in algebraic geometry. For compact K\"ahler manifolds, which form a natural transcendental analog of the projective world and where much of the minimal model program has recently been established, the question remains open in general.

	Let $\pi\colon\cX\to\Delta$ be a \emph{smooth family} of compact $n$-dimensional complex manifolds over the unit disk $\Delta\subset\C$ (i.e., $\pi$ is a proper surjective holomorphic submersion from a complex manifold) with fibers $X_t:=\pi^{-1}(t)$. Define the \emph{$m$-genus} of $X_t$ to be 
	\[
	P_m(X_t)=\dim_{\mathbb{C}}H^0(X_t,K^{\otimes m}_{X_t}),\qquad \text{for any integer $m\ge 1$}.
	\]
For compact K\"ahler manifolds, Siu posed the famous deformation invariance conjecture of plurigenera: 
	\begin{conjecture}[{\cite[Conjecture~0.4]{Siu02}, \cite[Conjecture~2.1]{SiuICM02}}]\label{conj:siu}
		For a smooth family $\pi: \cX \to \Delta$ of compact K\"ahler manifolds and every integer $m\ge1$, $P_m(X_t)$ is independent of $t\in\Delta$.
	\end{conjecture}

For smooth projective families this was proved by Siu, first for fibers of general type \cite{Si98} and then in general \cite{Siu02}. Siu's analytic proof \cite{Siu02} proceeds by constructing a singular hermitian metric on a pluricanonical bundle $(m-1)K_{\cX/\Delta}$ as a limit of metrics defined by inductively extended sections, using Siu's global generation together with a version of the Ohsawa--Takegoshi extension theorem \cite{OT87}. Kawamata \cite{Kaw99} and Takayama \cite{Tk07} extended the result to families whose fibers have canonical singularities, and P\u{a}un \cite{Pau07} gave a simplified proof of Siu's theorem. S. Rao--I.-H. Tsai \cite{RT1,RT2} proved fiberwise Moishezon family case, using Takayama's theorem and bimeromorphic embedding devoloped there. More recently, M.-L. Li--Rao--K.
Wang--M. Wang \cite{Ra26} can prove deformation invariance of plurigenera for a smooth family with uncountably many Moishezon fibers. In \cref{conj:siu}, the K\"ahler condition is fiberwise: each fiber $X_t$ is K\"ahler, while the total space $\cX$ is not required to be K\"ahler. J. Feng--Y. Li--Rao \cite[Theorem~3.8]{FLR26} shows that \cref{conj:siu} is equivalent to many other statements, for instance that every pluricanonical form on a fiber extends to the nearby fibers, or that the sheaves $R^1\pi_*\cO_\cX(mK_{\cX/\Delta})$ are torsion free.

	\subsection{Main results}

	Throughout the paper, we assume that good divisorial Zariski decompositions exist on compact K\"ahler manifolds, as stated below.

\begin{assumption}\label{prop:goodzariski}
	If $X$ is a compact K\"ahler manifold, $B$ is an snc $\Q$-divisor with coefficients in $[0,1]$, and $K_X+B$ is pseudo-effective, then there is a bimeromorphic modification $\mu\colon Y\to X$ from a compact K\"ahler manifold such that
		\[
		\mu^*(K_X+B)\sim_\Q P+R,\quad P\ \text{semiample},\quad R=N(\mu^*(K_X+B)) .
		\]
	Here \emph{semiample} means that, for some sufficiently large and divisible integer $m>0$, the line bundle $\mathcal{O}_Y(mP)$ is globally generated.
\end{assumption}
	
	This is announced in \cite{OAILog} under the assumption of logarithmic Iitaka subadditivity for projective varieties; the precise statement is recalled in \cref{sec:prelim}. We show that a suitable minimal model program on a compact K\"ahler manifold terminates with a good minimal model under the good Zariski decomposition assumption. Applied to a fiber $X_{t_0}$ of a family $\cX\to\Delta$, this minimal model program extends to a sequence of $K_{\cX}$-negative bimeromorphic contractions of the family over a neighbourhood of $t_0$, and this allows us to prove the conjecture of Siu on the invariance of plurigenera.

	\begin{theorem}[\cref{thm:main}]\label{thm:intro-main}
		Under the good divisorial Zariski decomposition \cref{prop:goodzariski}, let $\pi\colon\cX\to\Delta$ be a smooth family of compact K\"ahler manifolds. Then for every integer $m\ge1$, $P_m(X_t)$ is independent of $t\in\Delta$.
	\end{theorem}

	Equivalently, \cref{conj:siu} holds under this assumption. No K\"ahler form on the total space is needed, and the fibers are not assumed to be projective or Moishezon. The main new ingredient is the following termination result.

	\begin{theorem}[\cref{thm:mmp}]\label{thm:intro-mmp}
Under the good divisorial Zariski decomposition \cref{prop:goodzariski}, let  $X$ be a compact K\"ahler manifold with $K_X$ pseudo-effective. Then there is a finite sequence of $K_X$-negative divisorial contractions and flips $X\dashrightarrow M$ such that $M$ is a terminal, strongly $\Q$-factorial compact K\"ahler variety with $K_M$ semiample.
	\end{theorem}

	For projective varieties, the corresponding statement, that an MMP with scaling of an ample divisor terminates as soon as a good minimal model exists is due to Demailly--Hacon--P\u{a}un \cite{DHP13}, Hacon--C. Xu \cite{HX13} and C. Lai \cite{Lai11}. In the K\"ahler setting, we prove an analogous result under the assumption that good divisorial Zariski decompositions exist as stated in \cref{prop:goodzariski}. This is exactly what the proof of \cref{thm:intro-main} needs, as the extension theorem in \cref{cor:through} applies to the steps of the MMP.

		\subsection{Strategy of the proofs}

	We first explain the proof of \cref{thm:intro-mmp}. Let $X$ be a compact K\"ahler manifold with $K_X$ pseudo-effective. We run a $K_X$-MMP with scaling of a K\"ahler class, which exists by Hacon--L. Xie \cite{HX26} (see \cref{prop:mmp-scaling}). The proof contains two stages:

	In the first stage the MMP contracts the negative part of the divisorial Zariski decomposition $N(K_X)$. The MMP contract precisely the negative part $N(K_X)$ after finitely many steps by \cref{lem:prefix}. So either the MMP terminates, or after finitely many steps it reaches a terminal strongly $\Q$-factorial compact K\"ahler variety $S$ with $K_S$ pseudo-effective and $N(K_S)=0$. From then on, every further step of the MMP is a flip, which is isomorphic to $S$ in codimension one by \cref{lem:flips}.

	In the second stage, we show that these flips terminate. This is where the good divisorial Zariski decomposition enters. Choose a K\"ahler class $\omega$ on $S$ such that $K_S+\omega$ is K\"ahler, and run the MMP $S=S_0\dashrightarrow S_1\dashrightarrow\cdots$ with scaling of $\omega$, with thresholds $\lambda_0\ge\lambda_1\ge\cdots$. Using \cref{prop:goodzariski}, we construct a model $M$, isomorphic to $S$ in codimension one, with $K_M$ semiample in \cref{lem:smallmodel}. A relative MMP over the base of the fibration defined by $K_M$ then produces a model $Y$, again isomorphic to $S$ in codimension one, and a number $\delta>0$ such that $K_Y+s\omega_Y$ is nef for every $s\in[0,\delta]$ by \cref{lem:interval}. If the MMP with scaling were infinite, we would have $\lambda_i>0$ and $\lambda_i\to0$, so there would be an index $i$ with $0<\lambda_i<\lambda_{i-1}\le\delta$. For $t=\lambda_i$ and $t=\lambda_{i-1}$, the classes $K_{S_i}+t\omega_i$ and $K_Y+t\omega_Y$ are both nef. Since $S_i$ and $Y$ are isomorphic in codimension one, the negativity lemma as in \cref{lem:twoparam}, shows that $K_{S_i}$ is nef, which contradicts $\lambda_i>0$. Therefore the program terminates with a terminal, strongly $\Q$-factorial compact K\"ahler variety $M$ with $K_M$ nef, and $K_M$ is semiample by \cref{prop:abundance}.

	We now explain how \cref{thm:intro-mmp} implies \cref{thm:intro-main}. Fix $m\ge1$ and $t_0\in\Delta$. It is enough to show that $P_m(X_t)$ is constant for $t$ near $t_0$. If $P_m(X_{t_0})=0$, this follows from upper semicontinuity, so we may assume $P_m(X_{t_0})>0$. Then $K_{X_{t_0}}$ is pseudo-effective, and \cref{thm:intro-mmp} gives a $K_{X_{t_0}}$-MMP $X_{t_0}\dashrightarrow M_{t_0}$ with finitely many steps, where $M_{t_0}$ is terminal, compact K\"ahler, and $K_{M_{t_0}}$ is semiample. Two results below finish the proof.

	The first one is the extension theorem of Hacon--Xie \cite[Theorem~7.7]{HX26}, which generalizes Hacon--Li--Rao \cite[Proposition~4.16]{HLR26} and goes back to Koll\'ar \cite[Theorem~2]{Kol21}: every step of an MMP of a compact K\"ahler fiber with canonical singularities extends to the family over a neighbourhood of $t_0$ by \cref{thm:deformation}. Since our program has finitely many steps, after shrinking the base we obtain a flat proper family $\mathcal{M}\to\Delta$ with central fiber $M_{t_0}$. Its fibers $M_t$ have canonical singularities and each $M_t$ is bimeromorphic to $X_t$. Plurigenera are bimeromorphic invariants of compact varieties with canonical singularities, so $P_m(X_t)=P_m(M_t)$ for every $t$ near $t_0$, including $t=t_0$ by \cref{cor:through}.

	The second ingredient is the invariance of plurigenera theorem
of Levine \cite{Levine83,Levine85} and Koll\'ar--Mori
\cite[Theorem~12.5.5]{KM92}. Given a flat proper family, the
plurigenera are locally constant near a central fiber
with semiample canonical divisor and terminal singularities. Koll\'ar observed that the
same proof applies when the central fiber belongs to Fujiki's
class $\mathcal C$
\cite[comments after Conjecture~17]{Kol21}.
For the reader's convenience, we include a complete proof
in \cref{thm:lkm}. Let us briefly sketch the idea of the proof. For a sufficiently divisible
integer $r>0$, a general section of $\omega_{X_{t_0}}^{[r]}$ extends to
every infinitesimal neighbourhood of the central fiber.
The theorem of formal function then gives an extension over a
neighbourhood of $t_0$, defining a cyclic cover
$\mathcal Y\to\mathcal M$.
Its central fiber $Y_{t_0}$ belongs to Fujiki's class
$\mathcal C$ and has canonical singularity.
The $E_1$-degeneration and base change imply that
$h^n(Y_t,\mathcal O_{Y_t})$ is locally constant, where
$n=\dim M_{t_0}$. For cyclic cover,
$h^n(Y_t,\mathcal O_{Y_t})=\sum_{j=1}^r P_j(M_t)$. Since cohomological dimension is upper semi-continuous, each $P_j(M_t)$ is locally constant. Since $\mathcal{X} \to \mathcal{M}$ is fiberwise bimeromorphic, 
\[
P_m(X_t)=P_m(M_t)=P_m(M_{t_0})=P_m(X_{t_0}).
\]
Thus $P_m$ is constant constant on $\Delta$.

The K\"ahler hypothesis on the fiber $X_{t_0}$ is used to run a K\"ahler MMP by \cref{thm:intro-mmp}, and to apply the extension theorem of Hacon--Xie in \cref{thm:deformation}. In contrast with the analytic approach, our argument uses no $L^2$ extension theorem and no positive metric on the total space.

\subsection{Recent progress on invariance of plurigenera}

	Deformation invariance of plurigenera is a central problem in the deformation theory of compact complex manifolds, see Rao--Tsai \cite[Section~1.3]{RT2} and Feng--Li--Rao \cite[Introduction]{FLR26} for surveys of the main results. Let us recall the developments that is closely related to this paper.

	In the last few months there has been rapid progress. F. Deng--J. Ning--Z. Wang--X. Zhou \cite{DNWZ26} proved the invariance of plurigenera for K\"ahler families whose fibers have nef canonical bundle, and Chen--Rao--Wang \cite{CRW26} treated smooth families with one fiber having nef canonical bundle, as well as essentially Moishezon families. Y. Li--P\u{a}un--Wang--Zhou \cite{LPWZ26} proved the invariance of all plurigenera for smooth families with K\"ahler total space, without nefness assumption. Their proof is based on the plurisubharmonic variation of singular K\"ahler--Einstein metrics in the perturbed classes $c_1(K_{X_t})+\varepsilon[\omega_t]$, and the Ohsawa--Takegoshi theorem. Hacon--Li--Rao \cite{HLR26} established local stability of pseudo-effectiveness of canonical divisors and local constancy of volumes of adjoint classes in K\"ahler families when central fiber is projective. Feng--Li--Rao \cite{FLR26} studied several aspects of the invariance of plurigenera for K\"ahler varieties with canonical singularities. Most recently, K. Liu--X. Wan \cite{LiuWan26} extended the theorem of Li--P\u{a}un--Wang--Zhou to pairs. 

	All of these results assume either some positivity of the canonical class or a K\"ahler form on the total space. The assumption that the total space is K\"ahler is stronger than the fiberwise K\"ahler hypothesis of \cref{conj:siu}. \Cref{thm:intro-main} works with the fiberwise hypothesis only, at the cost of assuming the good divisorial Zariski decomposition, and its proof replies on the minimal model program rather than on $L^2$ method.

		\subsection*{Organization}
	In Section 2, we collect the background: singularities and positivity, Boucksom's divisorial Zariski decomposition, the K\"ahler MMP with scaling, and the extension of an MMP to nearby fibers. We prove the termination \cref{thm:intro-mmp} in Section 3. In the last section, we prove the theorem of Levine and Koll\'ar--Mori for central fibers in Fujiki's class $\mathcal C$, and then the main \cref{thm:intro-main}.

	\textbf{Acknowledgements.}
	Rao would first like to express his sincere gratitude to Professor Kefeng
Liu, who led him to the field of deforamtion theory and suggested Siu’s conjecture as his thesis topic, for his constant encouragement and help.
Together with Professor C. D. Hacon, we proved the deformation invariance of plurigenera for K\"ahler
families of threefolds in \cite{HLR26}, while with Jiawei Feng, we studied that for K\"ahler families of
compact K\"ahler threefolds with canonical singularities and maximal Albanese dimension and also
a log analogue of Koll\'ar’s deformation theorem of general type in \cite{FLR26}. Along this process we
realize that Koll\'ar’s method of extending the steps of a minimal model program to a family from
\cite{Kol21} is the right way to Siu’s conjecture for the fiberwise K\"ahler family. The remaining difficulty
was to run such a program with termination on a compact Kähler manifold. This is where the
paper of OpenAI on log abundance \cite{OAILog} caught our attention: it establishes abundance and
the existence of good divisorial Zariski decompositions for compact K\"ahler manifolds, which is
exactly what we need. We are very grateful to Christopher and Jiawei for these collaborations, from
which this project grew.
	
	\textbf{AI statments.} The proof of termination of K\"ahler MMP was initially  suggested by ChatGPT, while we haven't understood it so well and thus further complete a  proof in a more classical  approach. To be more precise, inspired by \cite{HX26} and \cite{DHP13,HX13,Lai11}, we construct the K\"ahler MMP with scaling and reduce the termination problem to the existence of good Zariski decomposition claimed in \cite{OAILog}, by basic properties of Boucksom--Zariski divisorial decomposition \cite{Bou04}.
	

	\section{Preliminaries}\label{sec:prelim}
	
	This section collects the preliminary results used in the paper and is organized into four parts. We first recall the basic notions of positivity, singularities, and plurigenera. We then review the properties of Boucksom's divisorial Zariski decomposition \cite{Bou04}, which plays a key role in our arguments. The last two parts are devoted to the K\"ahler minimal model program of \cite{HX26} and to the extension of an MMP from the central fiber of a family to nearby fibers of \cite{HX26} and \cite{Kol21,HLR26}; these are the technical tools in the proofs of our main results.
	
	\subsection{Singularities, positivity, and plurigenera}\label{subsec:basic}

We first recall the singularities for the generalized pairs, and then collect the tools needed to compare Bott--Chern classes on birational models, and we end this subsection by recalling the birational invariance of plurigenera and the deformation stability of canonical singularities.

A \emph{compact K\"ahler variety} is a compact complex space equipped with a K\"ahler form, which is a closed real $(1,1)$-form that, locally, is of the form $dd^c\varphi$, where $\varphi$ is a smooth strictly plurisubharmonic function. A compact complex space belongs to \emph{Fujiki's class $\mathcal{C}$} if it is bimeromorphic to a compact K\"ahler manifold.

	For a normal compact K\"ahler variety $X$, let $H^{1,1}_\BC(X,\R)$ denote the \emph{Bott--Chern cohomology} of real closed $(1,1)$-forms with local potentials, and write $\{D\}$ for the Bott--Chern class of a $\Q$-Cartier divisor $D$. When no confusion arises, we occasionally omit the braces and use $D$ to denote this class. Fix a K\"ahler form $\omega$ on $X$. A class $\alpha\in H^{1,1}_\BC(X,\R)$ is \emph{nef} if, for every $\varepsilon>0$, it has a smooth representative $\alpha_\varepsilon\ge-\varepsilon\omega$. It is \emph{pseudo-effective} if it contains a closed positive $(1,1)$-current, and \emph{big} if it contains a \emph{K\"ahler current}, that is, a closed positive current $T$ satisfying $T\ge\varepsilon\omega$ for some $\varepsilon>0$. $\alpha$ is \emph{modified big} if there are a modification $\mu\colon Y\to X$ and a big class $\beta$ on $Y$ such that $\mu_*\beta=\alpha$.

\begin{lemma}[{\cite[Lemma~2.38]{DHP22}}]\label{lem:nefpullback}
		Let $\mu_X\colon W\to X$ be a projective bimeromorphic morphism of normal compact K\"ahler varieties. A class $\alpha\in H^{1,1}_\BC(X,\R)$ is nef if and only if $\mu_X^*\alpha$ is nef.
	\end{lemma}

	We recall the notion of a generalized K\"ahler pair (cf. \cite[Definition~2.1]{HLR26}). Let $X$ be a normal compact K\"ahler variety, let $B\ge0$ be an $\mathbb R$-divisor on $X$, and let $\boldsymbol{\beta}$ be a real closed b-$(1,1)$-current, where \emph{b-$(1,1)$-current} is a compatible collection of closed real $(1,1)$-currents $\beta_{X'}$, one on each proper bimeromorphic model $X'\to X$, satisfying $p_*\beta_{X_1}=\beta_{X_2}$ whenever $p\colon X_1\to X_2$ is a bimeromorphic morphism between models of $X$. We call $\bbeta_{X'}$ the trace of $\boldsymbol{\beta}$ on $X'$ and write $\bbeta_X$ for its trace on $X$. The triple $(X,B+\boldsymbol{\beta})$ is a \emph{generalized K\"ahler pair} if there is a resolution $\nu\colon X'\to X$ such that $\boldsymbol{\beta}$ descends to $X'$, the Bott--Chern class of $\bbeta_{X'}$ is nef, and the adjoint class descends to $X$: there are an $\mathbb R$-divisor $B'$ with simple normal crossing support and $\nu_*B'=B$, and a class $\gamma\in H^{1,1}_{\mathrm{BC}}(X,\mathbb R)$ such that
	\[
	[K_{X'}+B'+\bbeta_{X'}]=\nu^*\gamma.
	\]
	Choose a log resolution $\mu_X: Y\to X$ dominating $X'$, with the exceptional locus and the support of the induced boundary $B_Y$ having simple normal crossings. Let $\bbeta_Y$ be the trace of $\boldsymbol{\beta}$ on $Y$. Then
	\[
	K_Y+B_Y+\bbeta_Y=\mu_X^*(K_X+B+\bbeta_X),\qquad (\mu_X)_*B_Y=B.
	\]
	The \emph{generalized discrepancy} of a divisor $E$ over $X$ is $a(E,X,B+\boldsymbol{\beta}):=-\operatorname{mult}_E(B_Y)$. The pair is \emph{generalized terminal} if these discrepancies are positive for every exceptional divisor over $X$. The generalized pair $(X,B+\bbeta)$ admits \emph{generalized canonical singularities} if $a(E,X,B+\bbeta)\ge 0$ for every exceptional divisor $E$ over $X$. The generalized pair $(X,B+\bbeta)$ is \emph{generalized klt (gklt)} if $a(E,X,B+\bbeta)>-1$ for every divisor $E$ over $X$. The generalized pair $(X,B+\bbeta)$ is \emph{generalized lc (glc)} if $a(E,X,B+\bbeta)\ge -1$ for every divisor $E$ over $X$. These conditions are independent of the chosen log resolution. The pair is \emph{generalized dlt (gdlt)} if it is log smooth on an open subset $U\subset X$, every prime divisor $E\subset U$ has $-1\le a(E,X,B+\boldsymbol{\beta})\le0$, and every divisor over $X$ centered in $X\setminus U$ has discrepancy greater than $-1$ (cf. \cite[Definitions~2.7 and 2.9]{das2026mmpgeneralizedpairskahler}). A normal compact K\"ahler variety is \emph{strongly $\Q$-factorial} if every reflexive coherent sheaf of rank one has an invertible reflexive power.

	We next introduce two tools for working with Bott--Chern classes: Demailly's support theorem and a criterion for Bott--Chern classes coming as pulled back from a birational model.

\begin{theorem}[Deamilly's support theorem]\label{thm:2ndsupp}
	Let $X$ be a compact complex manifold, 	$D$ be a divisor on $X$ whose support is $\text{supp}(D) = \bigcup_{i=1}^r D_i$, and $\Theta$ be a real closed $(1,1)$-current of order $0$ on $X$ whose support is contained in $\operatorname{Supp}(D)$. Then $$\Theta = \sum_{i=1}^r c_i [D_i],$$ where $[D_i]$ denotes the integral current induced by divisor $D_i$, and $c_i \in \mathbb{R}$.
\end{theorem}
\begin{proof}
Let
\[
A=\bigcup_{i=1}^r\operatorname{Sing}(D_i)
  \,\cup\,\bigcup_{i\ne j}(D_i\cap D_j).
\]
Then $A$ is an analytic subset of codimension at least two
in $X$, and the hypersurfaces $D_i\setminus A$ are smooth,
connected, and pairwise disjoint.

By \cite[III, Theorem 2.13]{agbook}, a closed current of order
zero supported on a connected smooth complex hypersurface
is a constant multiple of its current of integration.
Consequently, for each $i$ there is a real constant $c_i$
such that, on $X\setminus A$,
\[
\Theta=\sum_{i=1}^r c_i[D_i].
\]
The constants are real because $\Theta$ is real.

Thus
\[
R:=\Theta-\sum_{i=1}^r c_i[D_i]
\]
is a closed $(1,1)$-current of order zero supported on $A$.
Since $\operatorname{codim}_X A\ge2$, the first support
theorem \cite[III, Theorem 2.10]{agbook} gives $R=0$.
This proves the assertion.
\end{proof}

	\begin{lemma}[{\cite[Lemma 3.3]{HP16}}]\label{lem:pullbackBott}
	Let $f: X \rightarrow Y$ be a proper bimeromorphic morphism with connected fibers between normal compact complex spaces with rational singularities. Assume that $X$ and $Y$ are in Fujiki's class $\mathcal{C}$. Then, the pull back $f^*: H_{\mathrm{BC}}^{1,1}(Y) \rightarrow H_{\mathrm{BC}}^{1,1}(X)$ is injective, such that $$\operatorname{Im}\left(f^*\right)=\left\{\alpha \in H_{\mathrm{BC}}^{1,1}(X) \mid \alpha \cdot C=0 \text{ for all curves } C \subset X \text{ such that } f_*(C)=\mathrm{pt}\right\}.$$
	\end{lemma}

To compare Bott--Chern classes whose pullbacks differ by exceptional divisors, we use the following form of the negativity lemma.

\begin{lemma}[{Negativity lemma, \cite[Lemma 2.8]{das2026transcendentalminimalmodelprogram}}]\label{lem:negativity}
		Let $\mu_Y\colon W\to Y$ be a proper contraction morphism of normal complex analytic varieties, where $W$ is K\"ahler and $Y$ is relatively compact. Let $E$ be a $\mu_Y$-exceptional $\R$-divisor such that $-E$ is $\mu_Y$-nef. Then $E\ge0$.
	\end{lemma}

The negativity lemma gives the following useful result

The negativity lemma gives the following comparison for affine
families of Bott--Chern classes.

\begin{lemma}\label{lem:twoparam}
Let $p\colon W\to X$ and $q\colon W\to Y$ be projective
bimeromorphic morphisms from a compact K\"ahler manifold $W$ to
normal compact K\"ahler varieties $X$ and $Y$. Let
$\alpha_X,\gamma_X\in H^{1,1}_\BC(X,\R)$ and
$\alpha_Y,\gamma_Y\in H^{1,1}_\BC(Y,\R)$ satisfy
\[
p^*\alpha_X-q^*\alpha_Y=E,
\qquad
p^*\gamma_X-q^*\gamma_Y=F,
\]
where $E$ and $F$ are $\R$-divisors exceptional over both $X$ and $Y$.
Suppose that $\alpha_X+t\gamma_X$ and $\alpha_Y+t\gamma_Y$ are nef
for $t=a,b$, where $a\ne b$. Then $E=F=0$. In particular,
$\alpha_X+t\gamma_X$ is nef if and only if
$\alpha_Y+t\gamma_Y$ is nef for every $t\in\R$.
\end{lemma}

\begin{proof}
For $t=a,b$, the equality
\[
p^*(\alpha_X+t\gamma_X)-q^*(\alpha_Y+t\gamma_Y)=E+tF
\]
shows that $E+tF$ is $q$-nef and $-(E+tF)$ is $p$-nef.
By \cref{lem:negativity}, we obtain $E+tF\le0$ and $E+tF\ge0$.
Thus $E+aF=E+bF=0$, and $a\ne b$ gives $E=F=0$.
The last assertion follows from \cref{lem:nefpullback}.
\end{proof}


	For a normal compact complex space $X$ with canonical singularities, let $\mathcal{O}_X(mK_X)$ denote the $m$-th reflexive power of the canonical sheaf $\omega_X$. We put \[
	P_m(X)=\dim H^0(X,\mathcal{O}_X(mK_X))
	\] for $m\ge1$. 
	\begin{lemma}\label{lem:birinv}
		If $X$ has canonical singularities and $\mu_X\colon\widetilde X\to X$ is a resolution, then, for every $m\ge1$,
		\[
		(\mu_X)_*\mathcal{O}_{\widetilde X}(mK_{\widetilde X}) = \mathcal{O}_X(mK_X).
		\]
		In particular, if $X$ and $Y$ are bimeromorphic normal compact complex varieties with canonical singularities, then $P_m (X)=P_m(Y)$ for every $m\geq 1$.
	\end{lemma}

	\begin{proof}
		The assertion is local on $X$, so we may assume that $X$ is Stein so that the canonical divisor $K_X$ exists by Cartan A theorem. Fix a canonical divisor $K_{\widetilde X}$ so that $K_X=(\mu_X)_*K_{\widetilde X}$. Since $X$ has canonical singularities we can write
		\[
		K_{\widetilde X}=\mu_X^*K_X+E,
		\]
		with $E$ a $\mu_X$-exceptional $\mathbb{Q}$-divisor, and $E\ge0$. For an open set $U\subseteq X$,
		\[
		H^0(U,\mathcal{O}_X(mK_X))=\{f : \operatorname{div}(f)+mK_X\ge0 \text{ on } U\},
		\]
		\[
		H^0(\mu_X^{-1}U,\mathcal{O}_{\widetilde X}(mK_{\widetilde X}))=\{f : \operatorname{div}(\mu_X^*f)+mK_{\widetilde X}\ge0 \text{ on } \mu_X^{-1}U\},
		\]
		where $f$ runs over meromorphic functions on $U$; these $\mu^*_Xf$ will run through all meromorphic functions on $\mu_X^{-1}U$, since $\mu_X$ is bimeromorphic. We prove $(\mu_X)_*\mathcal{O}_{\tilde{X}}(m K_{\tilde{X}}) = \mathcal{O}_X(mK_X)$ by showing equality of the local sections. If $\operatorname{div}(\mu_X^*f)+mK_{\widetilde X}\ge0$, push this divisor forward by $\mu_X$. Pushforward preserves effectivity, $(\mu_X)_*\operatorname{div}(\mu_X^*f)=\operatorname{div}(f)$, and $(\mu_X)_*K_{\widetilde X}=K_X$. Hence $\operatorname{div}(f)+mK_X\ge0$. On the other hand, if $D:=\operatorname{div}(f)+mK_X\ge0$, then $D$ is an effective $\mathbb{Q}$-Cartier divisor, so $\mu_X^*D\ge0$. Therefore
		\[
		\operatorname{div}(\mu_X^*f)+mK_{\widetilde X}=\mu_X^*D+mE\ge0 .
		\]
		
		For the consequence, given $X\dashrightarrow Y$ bimeromorphic, take a common resolution $p\colon W\to X$ and $q\colon W\to Y$, gives $\dim H^0(X,\cO(mK_X))= \dim H^0(W,\cO(mK_W))= \dim H^0(Y,\cO(mK_Y))$.
	\end{proof}

	\begin{theorem}[{\cite{Kaw99}}]\label{thm:candef}
		Let $f: \cY\to\Delta$ be flat, and suppose the central fiber $Y_0 = f^{-1}(0)$ has canonical singularities. Then, in a neighbourhood of $Y_0$, the total space and the fibers have canonical singularities.
	\end{theorem}
	
	\subsection{Boucksom's divisorial Zariski decomposition}\label{subsec:negpart}
	
	We now turn to the divisorial negative part of a pseudo-effective class. After recalling Boucksom's construction and its behavior under modifications. At the end of this subsection, we formulate the good divisorial Zariski decomposition assumption and deduce its consequence.

Let $\alpha$ be a pseudo-effective class on a compact K\"ahler manifold $Y$ and $D\subset Y$ a prime divisor. If $\alpha$ is big, define
	\[
	\nu_D(\alpha):=\inf\{\nu(T,D)\mid T\in\alpha\text{ is a closed positive }(1,1)\text{-current}\},
	\]
	where $\nu(T,D)$ is the generic Lelong number of $T$ along $D$, that is $\nu(T,D):=\nu(T,x)$ for a general smooth point $x\in D$, where $\nu(T,x)$ is the usual Lelong number of $T$ at point $x$. If $\alpha$ is merely pseudo-effective, set \[
	\nu_D(\alpha)=\lim_{\varepsilon\downarrow0}\nu_D(\alpha+\varepsilon\omega),
	\] where $\omega$ is a K\"ahler class on $Y$; see \cite[Definition~3.1]{Bou04}. The \emph{negative part} of Boucksom's divisorial Zariski decomposition of $\alpha$ is
	\[
	N(\alpha)=\sum_D\nu_D(\alpha)\,D,
	\]
	a finite sum. 
	
	The following lemma collects the properties needed below. (i)--(iv) concern operations on a fixed manifold, (v)--(vii) describe the behavior of negative part under bimeromorphic modifications.
	
	\begin{lemma}\label{lem:zariski}
		Let $\alpha$ be a pseudo-effective class on a compact K\"ahler manifold $Y$.
		\begin{enumerate}[label=(\roman*)]
			\item\label{zar:calc} For every pseudo-effective class $\alpha'$, we have $N(\alpha+\alpha')\le N(\alpha)+N(\alpha')$ and $N(t\alpha)=tN(\alpha)$ for $t>0$. Moreover, $N(D)\le D$ for every effective $\R$-divisor $D$. The negative part $N(\alpha)$ has only finitely many components.
			\item\label{zar:nef} If $\beta$ is nef, then $N(\alpha+\beta)\le N(\alpha)$. In particular $N(\beta)=0$.
			\item\label{zar:neflimit}\label{zar:limit} If $\gamma$ is nef, then
			\[
			\lim_{t\downarrow0}\nu_D(\alpha+t\gamma)=\nu_D(\alpha)
			\]
			for every prime divisor $D$.
			\item\label{zar:subtract} If $F$ is an $\R$-divisor with $0\le F\le N(\alpha)$, then $\alpha-F$ is pseudo-effective and
			\[
			N(\alpha-F)=N(\alpha)-F.
			\]
			In particular $N(\alpha-N(\alpha))=0$.
			\item\label{zar:strict} If $\mu_Y\colon Y'\to Y$ is a modification from a compact K\"ahler manifold and $D'$ is the strict transform of a prime divisor $D\subset Y$, then $\nu_{D'}(\mu_Y^*\alpha)=\nu_D(\alpha)$ and $(\mu_Y)_*N(\mu_Y^*\alpha)=N(\alpha)$.
			\item\label{zar:exc} Let $\mu_X\colon Y\to X$ be a modification onto a normal compact K\"ahler variety, $\theta\in H^{1,1}_\BC(X,\R)$, and $F\ge0$ a $\mu_X$-exceptional $\R$-divisor such that $\mu_X^*\theta+F$ is pseudo-effective. Then $\mu_X^*\theta$ is pseudo-effective and
			\[
			N(\mu_X^*\theta+F)=N(\mu_X^*\theta)+F.
			\]
			\item\label{zar:pullup} If $\mu_Y\colon Y'\to Y$ is a modification from a compact K\"ahler manifold and $\alpha-N(\alpha)$ is nef, then $N(\mu_Y^*\alpha)=\mu_Y^*N(\alpha)$.
		\end{enumerate}
	\end{lemma}
	
	\begin{proof}
		(i) follows from \cite[Proposition~3.9 and Theorem~3.12]{Bou04} and the estimate in \cite[Proposition~3.6]{Bou04}:
		\[
			\nu_D(D_0)\le\nu([D_0],D)=\mult_DD_0.
		\]
		
		\smallskip
		(ii) A nef class has zero minimal multiplicity at every point \cite[Proposition~3.2]{Bou04}, so $N(\beta)=0$, and (i) gives $N(\alpha+\beta)\le N(\alpha)+N(\beta)=N(\alpha)$.
		
		\smallskip
		For (iii), by (ii) $\nu_D(\alpha+t\gamma)\le\nu_D(\alpha)$ for every $t>0$. The reverse inequality in the limit is the lower semicontinuity of $\nu_D$ on the pseudo-effective cone \cite[Proposition 3.5]{Bou04}.
		
		\smallskip
		(iv) Put Zariski projection of $\alpha$ as $Z(\alpha):=\alpha-N(\alpha)$, which is pseudo-effective with $N(Z(\alpha))=0$ \cite[Propositions~3.2 and~3.8]{Bou04}. Since
		\[
		\alpha-F=Z(\alpha)+N(\alpha)-F
		\]
		and $N(\alpha)-F\ge0$, the class $\alpha-F$ is pseudo-effective, and (i) gives
		\[
		N(\alpha-F)\le N(Z(\alpha))+N(N(\alpha)-F)\le N(\alpha)-F .
		\]
		Conversely, (i) gives $N(\alpha)\le N(\alpha-F)+N(F)\le N(\alpha-F)+F$.
		
		\smallskip
		(v) is \cite[Lemma~5.8]{Bou04} (see also \cite[Lemma~A.7]{das2026mmpgeneralizedpairskahler}).
		
		\smallskip
		(vi) is \cite[Lemma~A.5]{das2026mmpgeneralizedpairskahler}.
		
		\smallskip
		(vii) Put $\beta=\mu_Y^*\alpha$ and $Z(\alpha)=\alpha-N(\alpha)$. Since $\mu_Y^*Z(\alpha)$ is nef, (i) and (ii) give
		\[
		N(\beta)\le N(\mu_Y^*Z(\alpha))+N(\mu_Y^*N(\alpha))\le\mu_Y^*N(\alpha).
		\]
		If $D'$ is the strict transform of a prime divisor $D\subset Y$, then by (v)
		\[
		\mult_{D'}N(\beta)= \nu_{D'} (\beta)\stackrel{(v)}=\nu_D(\alpha)=\mult_DN(\alpha)=\mult_{D'}\mu_Y^*N(\alpha).
		\]
		Hence $G:=\mu_Y^*N(\alpha)-N(\beta)$ is an effective $\mu_Y$-exceptional divisor, and
		\[
		\beta-N(\beta)=\mu_Y^*Z(\alpha)+G.
		\]
		The left side has zero negative part by (iv). By (vi) applied to $\theta=Z(\alpha)$, the right side has negative part $N(\mu_Y^*Z(\alpha))+G=G$. Therefore $G=0$, that is, $N(\beta)=\mu_Y^*N(\alpha)$.
	\end{proof}
	
	\cref{lem:zariski} also allows us to define the negative part of Boucksom's divisorial Zariski decomposition on a singular space. If $X$ is a normal compact K\"ahler variety and $\alpha\in H^{1,1}_\BC(X,\R)$ is pseudo-effective, choose a resolution $\mu_X\colon Y\to X$ and set
\[
N(\alpha):=(\mu_X)_*N(\mu_X^*\alpha).
\]
By \cref{lem:zariski}.\ref{zar:strict}, this definition is independent of the resolution; see \cite[Definition~A.6 and Lemma A.7]{das2026mmpgeneralizedpairskahler}.

Under \cref{prop:goodzariski}, we can prove the abundance for klt pairs. 

\begin{theorem}\label{prop:abundance}
Under \cref{prop:goodzariski}, if $(X,B)$ is a klt pair, $X$ is a compact K\"ahler variety and $K_X+B$ is nef, then $K_X+B$ is semiample.
\end{theorem}

\begin{proof}
Let $f\colon X'\to X$ be a log resolution. Since $(X,B)$ is klt, we can write
	\[
K_{X'}+B'=f^*(K_X+B)+E,
	\]
where $B',E\geq0$ have no common components. $E$ is $f$-exceptional, and $B'$ is an snc divisor with coefficients in $[0,1)$. In particular, $K_{X'}+B'$ is pseudo-effective.

By \cref{prop:goodzariski}, there is a modification $\mu\colon Y\to X'$ from a compact K\"ahler manifold such that
	\[
\mu^*(K_{X'}+B')\sim_{\mathbb{Q}}P+R,\qquad P\ \text{semiample},\quad R=N\big(\mu^*(K_{X'}+B')\big).
	\]
Put $g:=f\circ\mu$. Then
	\[
\mu^*(K_{X'}+B')=g^*(K_X+B)+\mu^*E,
	\]
where $g^*(K_X+B)$ is nef and $\mu^*E\geq 0$ is $g$-exceptional, since $E$ is $f$-exceptional. By \cref{lem:zariski}, $R=\mu^*E$ as divisors. Hence
	\[
g^*(K_X+B)\sim_{\mathbb{Q}}P
	\]
is semiample. Therefore $K_X+B$ is semiample.
\end{proof}

\begin{remark}
In \cite[Theorem~2.13 and Proposition~2.7]{OAILog}, it is claimed that \cref{prop:goodzariski} holds if the logarithmic Iitaka subadditivity conjecture is true for projective varieties. A proof of Campana's orbifold Iitaka conjecture for rational snc boundaries on compact manifolds in Fujiki's class $\mathcal{C}$ has recently been announced in \cite[Theorem~1.1]{OAIOrbIitaka}. Logarithmic subadditivity is a special case. 

In the case of K\"ahler threefolds, log abundance has been established by Campana--H\"oring--Peternell and Das--Ou \cite{CHP16,DOI,DOII}.
\end{remark}

\subsection{The K\"ahler minimal model program}\label{subsec:kahlermmp}
For compact K\"ahler threefolds the minimal model program was established by H\"oring--Peternell \cite{HP16} and Das--Hacon \cite{das2024minimalmodelprogramkahler,das2024logminimalmodelprogram}, and extended to generalized pairs by Das--Hacon--Y\'a\~nez \cite{das2026mmpgeneralizedpairskahler}. In higher dimensions the transcendental K\"ahler MMP is developed by Hacon--Xie \cite{HX26}. We recall the parts of this theory that we use: the cone theorem, the contraction theorem, and the MMP with scaling of a K\"ahler class.

Let $X$ be a normal compact K\"ahler variety. Set $N^1_{\mathrm{an}}(X):=H^{1,1}_{\BC}(X,\R)$, and let $N_{1,\mathrm{an}}(X)$ be the space of real closed currents of bidimension $(1,1)$ on $X$, where two currents are identified if they have the same pairing with every class in $N^1_{\mathrm{an}}(X)$. In the projective setting, the extremal rays contracted by the MMP lie in the Mori cone $\overline{\mathrm{NE}}(X)$, the closed cone generated by the classes of curves (cf.\ \cite[Chapter~3]{KM98}). On a general compact K\"ahler variety a class in $N^1_{\mathrm{an}}(X)$ can be positive on every curve without being nef. One therefore replaces the Mori cone by the \emph{generalized Mori cone}
\[
\overline{\mathrm{NA}}(X):=\overline{\{[T]\in N_{1,\mathrm{an}}(X): T\text{ is a positive closed current of bidimension }(1,1)\}}.
\]

We use the following form of the cone theorem, due to Hacon--P\u{a}un \cite[Theorem~0.5]{HP24} (see also \cite[Corollary~3.5]{HLX26} and \cite[Theorem~1.3]{HX26}).

\begin{theorem}[{Cone theorem, \cite[Theorem 0.5]{HP24}}]\label{thm:kahler-cone}
Let $(X,B+\boldsymbol{\beta})$ be a compact K\"ahler generalized klt pair. Then there exists at most a countable collection of rational curves $\{C_i\}_{i\in I}$ on $X$ such that
$$\overline{\mathrm{NA}}(X)=\overline{\mathrm{NA}}(X)_{K_X+B+\bbeta_X\ge0}+\sum_{i\in I}\R_{\ge0}[C_i],$$
where $\overline{\mathrm{NA}}(X)_{K_X+B+\bbeta_X\ge0}$ is the subcone on which $K_X+B+\bbeta_X$ is nonnegative. Each $\R_{\ge0}[C_i]$ is a $(K_X+B+\bbeta_X)$-negative extremal ray, and
\[
0<-(K_X+B+\bbeta_X)\cdot C_i\le2\dim X.
\]
If $B+\bbeta_X$ is big, then $I$ is finite.
\end{theorem}

Let $(X,B+\boldsymbol{\beta})$ be a strongly $\Q$-factorial compact K\"ahler generalized klt pair. A projective bimeromorphic morphism $\phi\colon X\to Z$ with connected fibers onto a normal compact K\"ahler space is a \emph{$(K_X+B+\bbeta_X)$-negative contraction} if $-(K_X+B+\bbeta_X)$ is $\phi$-K\"ahler. It is \emph{extremal} if the curves contracted by $\phi$ are exactly those whose classes lie in a single extremal ray of $\overline{\mathrm{NA}}(X)$. An extremal negative contraction is \emph{divisorial} if its exceptional locus is a prime divisor, and \emph{small} if its exceptional locus has codimension at least two. The \emph{flip} of a small $(K_X+B+\bbeta_X)$-negative contraction $\phi$ is a small projective bimeromorphic morphism $\phi^+\colon X^+\to Z$ from a normal space such that $K_{X^+}+B^++\bbeta_{X^+}$ is $\phi^+$-K\"ahler, where $B^+$ is the push-forward of $B$.

The next result provides the contractions and flips needed to run the K\"ahler minimal model program. 
\begin{proposition}[{\cite[Proposition 4.3]{HX26}}]\label{prop:HXcont}
Let $(X,B+\boldsymbol{\beta})$ be a generalized klt pair on an $n$-dimensional strongly $\Q$-factorial compact K\"ahler variety $X$. Let $R=\mathbb{R}^{+}[\Sigma] \subset \overline{\operatorname{NA}}(X)$ be a $\left(K_X+B+\boldsymbol{\beta}_X\right)$-negative extremal ray. Then the contraction $f: X \rightarrow Y$ associated to $R$ exists and\\
(1) $f$ is projective and $Y$ is Kähler with rational singularities.\\
(2) If $f$ is a divisorial contraction, then $Y$ is strongly $\mathbb{Q}$-factorial. There is a canonical pushforward map:
$$
f_*: H_{\mathrm{BC}}^{1,1}(X) \rightarrow H_{\mathrm{BC}}^{1,1}(Y) \text { such that } f_* \circ f^*=\mathrm{Id}
$$
(3) If $f$ is a flipping contraction, then the flip $f^{+}: X^{+} \rightarrow Y$ exists, $X^{+}$ is strongly $\mathbb{Q}$-factorial, and there is a canonical isomorphism
$$
f_*: H_{\mathrm{BC}}^{1,1}(X) \rightarrow H_{\mathrm{BC}}^{1,1}\left(X^{+}\right) \text {induced by } \varphi: X \dashrightarrow X^{+}
$$
\end{proposition}

\begin{definition}\label{def:mmp}
Let $(X,B+\boldsymbol{\beta})$ be a strongly $\Q$-factorial compact
K\"ahler generalized klt pair. A \emph{$(K_X+B+\bbeta_X)$-minimal
model program} (or \emph{$(K_X+B+\bbeta_X)$-MMP}) is a finite or
infinite sequence of birational maps
$$
X=X_0\dashrightarrow X_1\dashrightarrow X_2\dashrightarrow\cdots,
$$
where $B_i,\bbeta_{X_i}$ are the pushforward of $B$ and $\bbeta_{X}$ onto $X_i$, such that each step is either a divisorial contraction of a $(K_{X_i}+B_i+\bbeta_{X_i})$-negative extremal ray, or  the flip of a flipping contraction of a $(K_{X_i}+B_i+\bbeta_{X_i})$-negative extremal ray.
\end{definition}

The desired outcome of the MMP when $K_X+B+\bbeta_X$ is pseudo-effective is a minimal model. We recall the precise definition.
\begin{definition}\label{def:gmm}
Let $(X,B+\boldsymbol{\beta})$ be a compact K\"ahler generalized klt pair. A \emph{minimal model} of $(X,B+\boldsymbol{\beta})$ is a bimeromorphic map $\varphi\colon X\dashrightarrow X'$ to a compact K\"ahler generalized klt pair $(X',B'+\boldsymbol{\beta})$ such that:
\begin{enumerate}[label=(\alph*)]
	\item $\varphi$ extracts no divisor, and $B'$ is the push-forward of $B$;
	\item $K_{X'}+B'+\bbeta_{X'}$ is nef;
	\item for a common resolution $\mu_X\colon W\to X$, $\mu_{X'}\colon W\to X'$, we have
	\[
	\mu_X^*(K_X+B+\bbeta_X)=\mu_{X'}^*(K_{X'}+B'+\bbeta_{X'})+F,
	\]
	where $F\ge0$ is $\mu_{X'}$-exceptional and its support contains the strict transform of every prime divisor contracted by $\varphi$.
\end{enumerate}
A minimal model is \emph{good} if $K_{X'}+B'+\bbeta_{X'}$ is semiample, that is, if the contraction associated with it exists.
\end{definition}

We now recall the MMP with scaling. Choose a K\"ahler class
$\omega_0$ on $X_0:=X$ such that $K_X+B+\bbeta_X+\omega_0$ is
K\"ahler. At the $i$-th stage, write
$\alpha_i:=K_{X_i}+B_i+\bbeta_{X_i}$, and let $\omega_i$ be the
push-forward of $\omega_0$ under the preceding steps. Define
the \emph{nef threshold} by
$$
\lambda_i:=\inf\{t\geq 0\mid \alpha_i+t\omega_i \text{ is nef}\}.
$$
If $\alpha_i$ is nef, the program stops. Otherwise, we choose a
$\alpha_i$-negative extremal ray $R_i\subset\overline{\mathrm{NA}}(X_i)$
such that
$$
(\alpha_i+\lambda_i\omega_i)\cdot R_i=0,
$$
and take its contraction. If the contraction is divisorial,
we let $X_{i+1}$ be its target; if it is flipping contraction, we let $X_{i+1}$
be its flip. We then define $\omega_{i+1}$ using the corresponding
push-forward map in \cref{prop:HXcont} and repeat the construction.
If the contraction is of fiber type, the program terminates
with a Mori fiber space.

The following proposition records this construction.

\begin{proposition}[{MMP with scaling, \cite{HX26}}]\label{prop:mmp-scaling}
Let $(X,B+\boldsymbol{\beta})$ be a strongly $\Q$-factorial compact K\"ahler generalized klt pair such that $K_X+B+\bbeta_X$ is pseudo-effective, and let $\omega$ be a K\"ahler class such that $K_X+B+\bbeta_X+\omega$ is K\"ahler. Then there is a $(K_X+B+\bbeta_X)$-MMP
\[
(X,B+\boldsymbol{\beta})=(X_0,B_0+\boldsymbol{\beta}_0)\dashrightarrow(X_1,B_1+\boldsymbol{\beta}_1)\dashrightarrow\cdots
\]
as in \cref{def:mmp} with the following properties. Here $\omega_i\in H^{1,1}_\BC(X_i,\R)$ denotes the push-forward of $\omega$, and the scaling thresholds are
\[
\lambda_{-1}:=1,\qquad \lambda_i:=\inf\{t\ge 0 \mid K_{X_i}+B_i+\bbeta_{X_i}+t\omega_i\text{ is nef}\}\quad(i\ge0).
\]
\begin{enumerate}[label=(\arabic*)]
	\item\label{prop:mmp-scaling:1} The thresholds satisfy $1=\lambda_{-1}\ge\lambda_0\ge\lambda_1\ge\cdots\ge0$. For $t\in[0,\lambda_{i-1}]$, the class $K_{X_i}+B_i+\bbeta_{X_i}+t\omega_i$ is nef if and only if $t\ge\lambda_i$.
	\item\label{prop:mmp-scaling:2} Every $(X_i,B_i+\boldsymbol{\beta}_i)$ is a strongly $\Q$-factorial compact K\"ahler generalized klt pair. The sequence terminates at $X_N$ if and only if $\lambda_N=0$. If $\lambda_i>0$, then $X_i\dashrightarrow X_{i+1}$ is the divisorial contraction or the flip of a $(K_{X_i}+B_i+\bbeta_{X_i})$-negative extremal ray $R_i$ with
	\[
	(K_{X_i}+B_i+\bbeta_{X_i}+\lambda_i\omega_i)\cdot R_i=0.
	\]
	\item\label{prop:mmp-scaling:3} For $t\in[0,\lambda_{i-1}]$ and a common resolution $\mu_X\colon W\to X$, $\mu_{X_i}\colon W\to X_i$, we have
	\begingroup
	\renewcommand{\theHequation}{scaledcomparison}
	\begin{equation}\tag{*}\label{eq:scaledcomparison}
		\mu_X^*(K_X+B+\bbeta_X+t\omega)=\mu_{X_i}^*(K_{X_i}+B_i+\bbeta_{X_i}+t\omega_i)+E_i(t),
	\end{equation}
	\endgroup
	where $E_i(t)\ge0$ is a $\mu_{X_i}$-exceptional divisor. Moreover, for $t\in[\lambda_i,\lambda_{i-1}]$,
	\[
	N\bigl(\mu_X^*(K_X+B+\bbeta_X+t\omega)\bigr)=E_i(t).
	\]
	\item\label{prop:mmp-scaling:4} For every $\varepsilon>0$, only finitely many $\lambda_i$ satisfy $\lambda_i\ge\varepsilon$. Hence either the sequence terminates with $K_{X_N}+B_N+\bbeta_{X_N}$ nef, or it is infinite and $\lambda_i\to0$.
\end{enumerate}
\end{proposition}

\begin{proof}
We construct the sequence by induction on $i$ and prove \hyperref[prop:mmp-scaling:1]{(1)} in \hyperlink{prop:mmp-scaling:step1}{Step 1}, \hyperref[prop:mmp-scaling:2]{(2)} in \hyperlink{prop:mmp-scaling:step2}{Step 2}, and \eqref{eq:scaledcomparison} in \hyperlink{prop:mmp-scaling:step3}{Step 3}. \hyperlink{prop:mmp-scaling:step4}{Step 4} proves the remaining assertion of \hyperref[prop:mmp-scaling:3]{(3)}, and \hyperlink{prop:mmp-scaling:step5}{Step 5} proves \hyperref[prop:mmp-scaling:4]{(4)}.

Let $\bar\omega$ denote the b-$(1,1)$ class that descends to $X$ with trace $\omega$. It is nef, so for every $\varepsilon\ge0$ the b-class $\boldsymbol{\beta}+\varepsilon\bar\omega$ is nef, and $(X_i,B_i+\boldsymbol{\beta}+\varepsilon\bar\omega)$ is a generalized pair with adjoint class $K_{X_i}+B_i+\bbeta_{X_i}+\varepsilon\omega_i$. 

The induction hypotheses are as follows: $X_i$, $B_i$ and $\omega_i$ have been constructed, and
\begin{enumerate}[label=(\alph*)]
	\item\label{prop:mmp-scaling:proof:a} $X_i$ is a strongly $\Q$-factorial compact K\"ahler variety and $X\dashrightarrow X_i$ is a bimeromorphic contraction;
	\item\label{prop:mmp-scaling:proof:b} $K_{X_i}+B_i+\bbeta_{X_i}+\lambda_{i-1}\omega_i$ is nef;
	\item\label{prop:mmp-scaling:proof:c} \eqref{eq:scaledcomparison} holds for $t\in[0,\lambda_{i-1}]$.
\end{enumerate}
For $i=0$ they hold with $\lambda_{-1}=1$ and $E_0=0$.

\hypertarget{prop:mmp-scaling:step1}{Step 1}. We prove \hyperref[prop:mmp-scaling:1]{(1)} for index $i$. By definition of $\lambda_i$ and hypothesis \hyperref[prop:mmp-scaling:proof:b]{(b)}, the class $K_{X_i}+B_i+\bbeta_{X_i}+t\omega_i$ is nef for $t\in[\lambda_i,\lambda_{i-1}]$, and it is not nef for $t\in[0,\lambda_i)$ and $\lambda_i \le \lambda_{i-1}$.

\hypertarget{prop:mmp-scaling:step2}{Step 2}. We next prove \hyperref[prop:mmp-scaling:2]{(2)} for the index $i$, and \hyperref[prop:mmp-scaling:proof:a]{(a)} and \hyperref[prop:mmp-scaling:proof:b]{(b)} for $X_{i+1}$.

If $\lambda_i=0$, then $K_{X_i}+B_i+\bbeta_{X_i}$ is nef and the sequence terminates at $X_i$. Assume from now on that $\lambda_i>0$, so that $K_{X_i}+B_i+\bbeta_{X_i}$ is not nef by \hyperlink{prop:mmp-scaling:step1}{Step 1}, and fix $\varepsilon\in(0,\lambda_i)$. The class $\omega_i$ is big, being the push-forward of a K\"ahler class under the bimeromorphic contraction $X\dashrightarrow X_i$, and $B_i+\bbeta_{X_i}$ is pseudo-effective. So $(X_i,B_i+\boldsymbol{\beta}+\varepsilon\bar\omega)$ is a generalized klt pair with big boundary $B_i+\bbeta_{X_i}+\varepsilon\omega_i$. Since
\[
	K_{X_i}+B_i+\bbeta_{X_i}+\varepsilon\omega_i+s(1-\varepsilon)\omega_i=K_{X_i}+B_i+\bbeta_{X_i}+\bigl(\varepsilon+s(1-\varepsilon)\bigr)\omega_i ,
\]
so that the nef threshold of this new pair with respect to $(1-\varepsilon)\omega_i$ is $$s_i = \inf\{\tau \ge 0 \mid K_{X_i}+B_i + \bbeta_i + \varepsilon \omega_i + \tau (1- \varepsilon)\omega_i \text{ is nef}\} = (\lambda_i-\varepsilon)/(1-\varepsilon)>0$$and that$$K_{X_i}+B_i + \bbeta_i + \varepsilon \omega_i + s_i (1- \varepsilon)\omega_i = K_{X_i}+B_i+\bbeta_{X_i}+\lambda_i\omega_i$$
Apply \cref{thm:kahler-cone} to $(X_i,B_i+\boldsymbol{\beta}_i+\varepsilon\omega_i)$. Since $B_i+\bbeta_{X_i}+\varepsilon\omega_i$ is big, the cone
theorem gives
\[
\overline{\mathrm{NA}}(X_i)
=
\overline{\mathrm{NA}}(X_i)_{
K_{X_i}+B_i+\bbeta_{X_i}+\varepsilon\omega_i\ge0}
+\sum_{j=1}^{r}\R_{\ge0}[C_j],
\]
where each $C_j$ spans a
$(K_{X_i}+B_i+\bbeta_{X_i}+\varepsilon\omega_i)$-negative
extremal ray.

For $t\in[\varepsilon,\lambda_{i-1}]$, the class
$K_{X_i}+B_i+\bbeta_{X_i}+t\omega_i$ is nonnegative on the
first summand. Indeed, it is a convex combination of the
class at $t=\varepsilon$ and the nef class at
$t=\lambda_{i-1}$. Therefore, $K_{X_i}+B_i+\bbeta_{X_i}+t\omega_i$ for $t\in[\varepsilon, \lambda_{i-1}]$ is nef if it has non-negative intersection with these finitely many curves $C_j$.

At $t=\lambda_i$, at least one of these intersections must
be zero. Otherwise, they would all remain positive for
$t<\lambda_i$ sufficiently close to $\lambda_i$.
Since $\lambda_i>\varepsilon$, the corresponding class
would then be nef, contrary to the infimum of $\lambda_i$.
Choose such a curve $\Sigma_i$ and set
$R_i=\R_{\ge0}[\Sigma_i]$. Then
\[
(K_{X_i}+B_i+\bbeta_{X_i}+\varepsilon\omega_i)\cdot\Sigma_i<0,
\qquad
(K_{X_i}+B_i+\bbeta_{X_i}+\lambda_i\omega_i)\cdot\Sigma_i=0.
\]
This implies that $(\lambda_i-\varepsilon)\,\omega_i\cdot\Sigma_i>0$. Hence
\begin{equation}\label{eq:signs}
	\omega_i\cdot\Sigma_i>0,\qquad (K_{X_i}+B_i+\bbeta_{X_i})\cdot\Sigma_i=-\lambda_i\,\omega_i\cdot\Sigma_i<0 .
\end{equation}
So $R_i$ is a $(K_{X_i}+B_i+\bbeta_{X_i})$-negative extremal ray. By \cref{prop:HXcont}, the contraction $f_i\colon X_i\to Z_i$ of $R_i$ exists, and its flip exists if the contraction is small. The contraction is bimeromorphic because the adjoint class of the perturbed pair is pseudo-effective, so it cannot be negative on the fibers of a fiber type contraction.

If $f_i$ is a divisorial contraction, put $X_{i+1}:=Z_i$ and $f_i^+:=\mathrm{id}$. If $f_i$ is a flipping contraction, let $f_i^+: X_{i+1}\to Z_i$ be its flip as a contraction of $(X_i,B_i+\boldsymbol{\beta}+\varepsilon\bar\omega)$. Let $\varphi_i: X_i\dashrightarrow X_{i+1}$ be the induced bimeromorphic map, and put $B_{i+1}:=\varphi_{i*}B_i$ and $\omega_{i+1}:=\varphi_{i*}\omega_i$. In both cases $X_{i+1}$ is strongly $\Q$-factorial and compact K\"ahler by \cref{prop:HXcont} and \cite[Corollary~1.8]{HX26}, and $X\dashrightarrow X_{i+1}$ is a bimeromorphic contraction. This is \hyperref[prop:mmp-scaling:proof:a]{(a)} for $X_{i+1}$.

By \cref{lem:pullbackBott}, we can write $K_{X_i}+B_i+\bbeta_{X_i}+\lambda_i\omega_i=f_i^*\delta_i$ for a class $\delta_i$ on $Z_i$, and $\delta_i$ is nef by \cref{lem:nefpullback}. Hence
\begin{equation}\label{id-delta}
K_{X_{i+1}}+B_{i+1}+\bbeta_{X_{i+1}}+\lambda_i\omega_{i+1}=f_i^{+*}\delta_i
\end{equation}
is nef. This is \hyperref[prop:mmp-scaling:proof:b]{(b)} for $X_{i+1}$.

It remains to see that $\varphi_i$ is a step of a $(K_{X_i}+B_i+\bbeta_{X_i})$-MMP, and in fact of the MMP of $(X_i,B_i+\boldsymbol{\beta}+\varepsilon'\bar\omega)$ with scaling of $(1-\varepsilon')\omega_i$ for every $\varepsilon'\in[0,\lambda_i)$. By \eqref{eq:signs}, the ray $R_i$ is negative for $K_{X_i}+B_i+\bbeta_{X_i}+\varepsilon'\omega_i$, and it is trivial for $K_{X_i}+B_i+\bbeta_{X_i}+\lambda_i\omega_i$. This settles the divisorial contraction case. In the flipping contraction case, the identity \eqref{id-delta} shows that over $Z_i$
\[
K_{X_{i+1}}+B_{i+1}+\bbeta_{X_{i+1}}+\varepsilon'\omega_{i+1}\equiv_{Z_i}\Bigl(1-\frac{\varepsilon'}{\lambda_i}\Bigr)(K_{X_{i+1}}+B_{i+1}+\bbeta_{X_{i+1}}).
\]
For $\varepsilon'=\varepsilon$ the left-hand side is $f_i^+$-K\"ahler by \cref{prop:HXcont}. So $K_{X_{i+1}}+B_{i+1}+\bbeta_{X_{i+1}}$ is $f_i^+$-K\"ahler, and then the left-hand side is $f_i^+$-K\"ahler for every $\varepsilon'\in[0,\lambda_i)$. Thus $f_i^+$ is the flip of $f_i$ for each of these pairs with $\varepsilon' \in [0,\lambda_i)$. This completes the proof of \hyperref[prop:mmp-scaling:2]{(2)} for the index $i$.

\hypertarget{prop:mmp-scaling:step3}{Step 3}. We prove \hyperref[prop:mmp-scaling:proof:c]{(c)} for $X_{i+1}$, that is, \eqref{eq:scaledcomparison} for $X_{i+1}$ and $t\in[0,\lambda_i]$. Let $W$ be a common resolution of $X$, $X_i$ and $X_{i+1}$. As the commutative diagram indicate below:
\begin{center}
	\begin{tikzcd}
	W & {X_{i+1}} & \\
	X & {X_i} & {Z_i}
	\arrow["{\mu_{X_{i+1}}}", from=1-1, to=1-2]
	\arrow["{\mu_X}"', from=1-1, to=2-1]
	\arrow["{\mu_{X_i}}"{description}, from=1-1, to=2-2]
	\arrow["{f_{i}^+}", from=1-2, to=2-3]
	\arrow[dashed, from=2-1, to=2-2]
	\arrow["{\varphi_i}"', dashed, from=2-2, to=1-2]
	\arrow["{f_i}"', from=2-2, to=2-3].
\end{tikzcd}
\end{center}

For $j=i,i+1$ write $K_W+\bbeta_W+\Delta_j=\mu_{X_j}^*(K_{X_j}+B_j+\bbeta_{X_j})$, so that $(\mu_{X_j})_*\Delta_j=B_j$, and put
\[
G_i:=\Delta_i-\Delta_{i+1}=\mu_{X_i}^*(K_{X_i}+B_i+\bbeta_{X_i})-\mu_{X_{i+1}}^*(K_{X_{i+1}}+B_{i+1}+\bbeta_{X_{i+1}}).
\]
We claim that $G_i$ is $\mu_{X_{i+1}}$-exceptional and effective. Since $(\mu_{X_{i+1}})_*\Delta_i=\varphi_{i*}B_i=B_{i+1}=(\mu_{X_{i+1}})_*\Delta_{i+1}$, so $G_i$ is exceptional over $X_{i+1}$, hence over $Z_i$. Moreover $-G_i$ is nef over $Z_i$, as the class $-(K_{X_i}+B_i+\bbeta_{X_i})$ is $f_i$-K\"ahler. Therefore $G_i\ge0$ by the negativity \cref{lem:negativity}.

Let $\gamma\in H^{1,1}_\BC(X_i,\R)$ and put $c:=(\gamma\cdot\Sigma_i)/((K_{X_i}+B_i+\bbeta_{X_i})\cdot\Sigma_i)$. Note that $(\gamma-c(K_{X_i}+B_i+\bbeta_{X_i})) \cdot R_i = 0$, so  $$ \gamma-c(K_{X_i}+B_i+\bbeta_{X_i})=f_i^*\delta$$for some $\delta\in H^{1,1}_\BC(Z_i,\R)$ by \cref{lem:pullbackBott}, and then $\varphi_{i*}\gamma=f_i^{+*}\delta+c(K_{X_{i+1}}+B_{i+1}+\bbeta_{X_{i+1}})$. Therefore
\begin{equation}\label{eq:gamma-comparison}
\mu_{X_i}^*\gamma-\mu_{X_{i+1}}^*\varphi_{i*}\gamma=c\,G_i .
\end{equation}
Now let $t\in[0,\lambda_i]$ and take $\gamma=K_{X_i}+B_i+\bbeta_{X_i}+t\omega_i$. Then $\varphi_{i*}\gamma=K_{X_{i+1}}+B_{i+1}+\bbeta_{X_{i+1}}+t\omega_{i+1}$, and $c=1-t/\lambda_i$ by \eqref{eq:signs}. By the induction hypothesis \hyperref[prop:mmp-scaling:proof:c]{(c)}, \eqref{eq:scaledcomparison} holds for $X_i$ at this value of $t$, because $t\le\lambda_i\le\lambda_{i-1}$. Combining \eqref{eq:scaledcomparison} for $X_i$ with \eqref{eq:gamma-comparison} gives \eqref{eq:scaledcomparison} for $X_{i+1}$, with
\[
	E_{i+1}(t)=E_i(t)+\Bigl(1-\frac{t}{\lambda_i}\Bigr)G_i .
\]
This divisor is effective because $E_i(t)\ge0$, $t\le\lambda_i$ and $G_i\ge0$. It is $\mu_{X_{i+1}}$-exceptional because $G_i$ is, and because $E_i(t)$ is $\mu_{X_i}$-exceptional.

\hypertarget{prop:mmp-scaling:step4}{Step 4}. We prove the second assertion of \hyperref[prop:mmp-scaling:3]{(3)}. Let $t\in[\lambda_i,\lambda_{i-1}]$, so that $K_{X_i}+B_i+\bbeta_{X_i}+t\omega_i$ is nef by \hyperref[prop:mmp-scaling:2]{(2)}. So \cref{lem:zariski} applies, and together with \eqref{eq:scaledcomparison} for $X_i$ it gives
\[
N\bigl(\mu_X^*(K_X+B+\bbeta_X+t\omega)\bigr)=N\bigl(\mu_{X_i}^*(K_{X_i}+B_i+\bbeta_{X_i}+t\omega_i)+E_i(t)\bigr)=E_i(t).
\]

\hypertarget{prop:mmp-scaling:step5}{Step 5}. We prove \hyperref[prop:mmp-scaling:4]{(4)} at the final step. Suppose, for a contradiction, that there is $\varepsilon>0$ such that $\lambda_i\ge\varepsilon$ for infinitely many $i$. Since the thresholds are non-increasing, the sequence is infinite and $\lambda_i\ge\varepsilon$ for all $i$. Fix $\varepsilon'\in(0,\varepsilon)$. The pair $(X,B+\boldsymbol{\beta}+\varepsilon'\bar\omega)$ is generalized klt with $B+\boldsymbol{\beta}+\varepsilon'\bar\omega$ big. The class $K_X+B+\bbeta_X+\varepsilon'\omega+(1-\varepsilon')\omega$ is K\"ahler. By the last part of \hyperlink{prop:mmp-scaling:step2}{Step 2}, each map $\varphi_i\colon X_i\dashrightarrow X_{i+1}$ is a step of the MMP of $(X,B+\boldsymbol{\beta}+\varepsilon'\bar\omega)$ with scaling of $(1-\varepsilon')\omega$. So this MMP with scaling is infinite, which contradicts to the termination result of \cite[Theorem~1.5]{HX26}. Therefore only finitely many $\lambda_i$ are at least $\varepsilon$. If the MMP is infinite, it follows that $\lambda_i\to0$.
\end{proof}

We conclude this subsection with an existence of MMP result. When the generalized boundary is big and the adjoint class is pseudo-effective over the base, the results of \cite{HX26} give a relative minimal model.
	\begin{proposition}\label{prop:relmodel}
		Let $(X,B+\beta)$ be a generalized klt pair on a strongly $\Q$-factorial compact K\"ahler variety $X$, such that $B+\bbeta_X$ is big. Let $g\colon X\to Z$ be a proper morphism to a normal compact K\"ahler variety and suppose $\alpha=K_X+B+\bbeta_X$ is pseudo-effective over $Z$. Then there is a bimeromorphic contraction $\varphi\colon X\dashrightarrow Y$ over $Z$, with $Y$ a strongly $\Q$-factorial compact K\"ahler variety, such that the following hold:\\
		\hypertarget{prop:relmodel:1}{$(1)$}
		\[
		\alpha_Y =K_Y+B_Y +\bbeta_Y = f_*(K_X+B+\beta)
		\]
		is the adjoint class of the generalized klt pair $(Y,B_Y+\beta)$,\\
		\hypertarget{prop:relmodel:2}{$(2)$} $\alpha_Y$ is nef over $Z$, and, \\
		\hypertarget{prop:relmodel:3}{$(3)$} on a common resolution $\mu_X\colon W'\to X$, $\mu_Y\colon W'\to Y$,
		\[
		\mu_X^*\alpha=\mu_Y^*\alpha_Y+E,\qquad E\ge0\ \text{exceptional over } Y,
		\]
		with $E$ positive along the strict transform of every prime divisor contracted by $\varphi$.
	\end{proposition}
	
	\begin{proof}
	By \cite[Theorem~5.14.(1)]{HX26}, the pair $(X,B+\bbeta)$ has a log terminal model $\varphi\colon X\dashrightarrow Y$ over $Z$ in the sense of \cite[Definition~2.15]{HX26}. In particular
	$$a(P,X,B+\beta)<a(P,Y,B_Y+\beta)\quad\text{for every prime divisor }P\subset X\text{ contracted by }\varphi.$$
	Thus $\varphi$ is a bimeromorphic contraction over $Z$, $Y$ is strongly $\Q$-factorial, $\alpha_Y=K_Y+B_Y+\bbeta_Y$ is nef over $Z$ with $B_Y=f_*B$, and $K_Y+B_Y+\bbeta_Y$ is gklt. 
	
It remains to prove (3). Choose a common resolution $\mu_X\colon W'\to X$, $\mu_Y\colon W'\to Y$ and define $B_{W'}$ and $B^Y_{W'}$ by
\[
K_{W'}+B_{W'}+\bbeta_{W'}=\mu_X^*\alpha,\qquad K_{W'}+B^Y_{W'}+\bbeta_{W'}=\mu_Y^*\alpha_Y,
\]
put $E=B_{W'}-B^Y_{W'}$. Then $\mu_X^*\alpha=\mu_Y^*\alpha_Y+E$, and
\[
\mult_PE=a(P,Y,B_Y+\beta)-a(P,X,B+\beta)\quad\text{for every prime divisor }P\subset W'.
\]
Since $\varphi^{-1}$ contracts no divisor, a prime divisor $P\subset W'$ that is not $\mu_Y$-exceptional is the strict transform of a divisor on both $X$ and $Y$, so $\mult_PE=0$ because $B_Y=f_*B$. Thus $E$ is $\mu_Y$-exceptional. The remaining prime divisors that are not $\mu_X$-exceptional are the strict transforms of the divisors contracted by $\varphi$, and along these $\mult_PE>0$ by the definition of a log terminal model. Hence $(\mu_X)_*E\ge0$. Moreover $-E$ is $\mu_X$-nef: a curve $C\subset W'$ contracted by $\mu_X$ maps to a point of $Z$, so $-E\cdot C=\alpha_Y\cdot(\mu_Y)_*C\ge0$ because $\alpha_Y$ is nef over $Z$. Therefore $E\ge0$ by the negativity \cref{lem:negativity}.\end{proof}

	\subsection{Extension of the minimal model program}\label{subsec:extension}
	
	The following result of Hacon--Xie \cite[Theorem~7.7]{HX26} lets us extend a K\"ahler MMP on the central fiber to the family. It generalizes \cite[Proposition~4.16]{HLR26}, where the central fiber is assumed projective. We state it in the form that boundary and nef part are zero, which is all we need.

	\begin{theorem}[a partial result of {\cite[Theorem 7.7]{HX26}, \cite[Proposition 4.16]{HLR26}}]\label{thm:deformation}
	Let $g\colon\mathcal{X}\to S$ be a flat proper morphism with connected fibers, where $\mathcal{X}$ is a normal complex variety and $S$ is a smooth, connected, relatively compact curve. Fix $0\in S$ and assume that the central fiber $X_0$ is a compact K\"ahler variety with canonical singularities. Then every sequence of $K_{X_0}$-MMP steps
\[
X_0\dashrightarrow X_0^{(1)}\dashrightarrow X_0^{(2)}\dashrightarrow\cdots
\]
extends over some open neighborhood $U\subset S$ of $0$ to a sequence of $K_{\mathcal{X}}$-negative proper bimeromorphic maps
\[
\mathcal{X}/U\dashrightarrow\mathcal{X}^{(1)}/U\dashrightarrow\mathcal{X}^{(2)}/U\dashrightarrow\cdots,
\]
where the fiber of $\mathcal{X}^{(i)}$ over $0$ is $X_0^{(i)}$ for every $i\ge0$.
	\end{theorem}
	
\begin{corollary}\label{cor:through}
		Let $\pi\colon\cX\to\Delta$ be a flat proper contraction morphism from a normal complex variety, and suppose that the central fiber $X_0$ is a compact K\"ahler variety with canonical singularities. Let
		\[
		X_0=X_0^0\dashrightarrow X_0^1\dashrightarrow\cdots\dashrightarrow X_0^N\dashrightarrow \cdots 
		\]
		be a sequence of steps of a $K_{X_0}$-MMP. Suppose that over a neighbourhood $U$ of $0$ this sequence extends to $K_{\mathcal{X}_U}$-negative bimeromorphic contractions
		\[
		\cX_U=\cX^0\dashrightarrow\cX^1\dashrightarrow\cdots\dashrightarrow\cX^N\dashrightarrow \cdots,
		\]
		where the fiber of $\cX^i$ over $0$ is $X_0^i$. Write $X_t^i$ for the fiber of $\cX^i$ over $t$. After shrinking $U$, fibers of $\mathcal{X}^i\to U$ have canonical singularities, and
		$$P_m(X_t)=P_m(X_t^i)\qquad\text{for all } m\ge1,\ 0\le i\le N,\ t\in U .$$
\end{corollary}
	
\begin{proof}
The central fiber $X_0^i$ of $\cX^i$ has canonical singularities because steps of a $K_{X_0}$-MMP preserve them \cite[Lemma~3.38]{KM98}. By \cref{thm:candef}, after shrinking $U$ every fiber $X_t^i$ has canonical singularities. By \cite[Lemma~2.22]{CRT25}, after further shrinking $U$, the induced maps $X_t^i\dashrightarrow X_t^{i+1}$ are bimeromorphic for $t$ near $0$, so \cref{lem:birinv} gives $P_m(X_t^i)=P_m(X_t^{i+1})$ for all $m\ge1$.
\end{proof}
	

	\section{Termination of the MMP with scaling}\label{sec:mmp}\label{subsec:negmmp}
	
	In this section we prove \cref{thm:intro-mmp}, restated below as \cref{thm:mmp}. Our argument is motivated by the results of Demailly--Hacon--P\u{a}un \cite[Theorem~2.4]{DHP13}, Hacon--Xu \cite[Corollary~2.9]{HX13}, and Lai \cite[Proposition~2.5]{Lai11}. Differently from these results, instead of assuming the existence of good minimal model, we assume only the good divisorial Zariski decomposition in \cref{thm:mmp}.
\begin{theorem}\label{thm:mmp}
Under the good divisorial Zariski decomposition \cref{prop:goodzariski}, let $X$ be a compact K\"ahler manifold with $K_X$ pseudo-effective. Then there is a finite sequence of $K_X$-negative divisorial contractions and flips $X\dashrightarrow M$ such that $M$ is a terminal, strongly $\Q$-factorial compact K\"ahler variety with $K_M$ semiample.
	\end{theorem}

	\begin{remark}
		In the projective setting, Hacon--Moraga \cite[Theorem~1]{HM20} prove that, assuming termination of flips for generalized klt pairs in dimensions at most $n-1$, every MMP for an $n$-dimensional generalized lc pair admitting a weak Zariski decomposition terminates. Their proof uses the ACC for generalized log canonical thresholds of \cite[Theorem~1.5]{BZ16}, a generalization of the ACC for log canonical thresholds established by Hacon--McKernan--Xu \cite{HMX14}. Their argument does not by itself provide the semiampleness needed for a good minimal model. Since our argument requires a good minimal model, we adopt the good divisorial Zariski decomposition \cref{prop:goodzariski}.
	\end{remark}

	Before proving the main result (\cref{thm:mmp}) of this section, let us first establish several preliminary results.
	\begin{lemma}\label{lem:contractedN}
		Let $X$ be a strongly $\Q$-factorial compact K\"ahler variety with klt singularities and $K_X$ pseudo-effective. Let $\omega$ be a K\"ahler class such that $K_X+\omega$ is K\"ahler, and let $X\dashrightarrow Y$ be a composition of finite steps of a $K_X$-MMP with scaling of $\omega$. Then $K_{Y}$ is pseudo-effective, and for every prime divisor $E\subset X$:
		\begin{enumerate}[label=(\roman*)]
			\item\label{lem:contractedN:contracted} If $E$ is contracted by $X\dashrightarrow Y$, then $\nu_E(K_X)>0$.
			\item\label{lem:contractedN:retained} If $E$ is not contracted and $E'\subset Y$ be the strict transform of $E$, then \[
			\nu_{E'}(K_{Y})=\nu_E(K_X).
			\]
		\end{enumerate}
	\end{lemma}
	
	\begin{proof}
		On a common resolution $\mu_X\colon W\to X$, $\mu_Y\colon W\to Y$, we have
		\[
		\mu_X^*K_X\sim_\Q \mu_Y^*K_{Y}+F,\qquad F\ge0\ \text{is $\mu_Y$-exceptional}.
		\]
		For a prime divisor $E\subset X$, let $\widetilde E\subset W$
be its strict transform. If $E$ is contracted by
$X\dashrightarrow Y$, then its discrepancy strictly increases,
so $\operatorname{mult}_{\widetilde E}F>0$.
If $E$ is not contracted by $X\dashrightarrow Y$, then it has
a strict transform $E'\subset Y$, and $\mu_Y$ maps
$\widetilde E$ birationally onto $E'$. Thus $\widetilde E$
is not $\mu_Y$-exceptional, and hence
$\operatorname{mult}_{\widetilde E}F=0$. By \cref{lem:zariski},
$$N(\mu_X^*K_X)=N(\mu_Y^*K_{Y})+F.$$
Let $\widetilde E\subset W$ be the strict transform of $E$. By \cref{lem:zariski}.\ref{zar:strict}, the coefficient of $\widetilde E$ in $N(\mu_X^*K_X)$ is $\nu_E(K_X)$, and we compare it with the coefficient of $\widetilde E$ on the right-hand side of the identity above. If $E$ is contracted by $X\dashrightarrow Y$, then $\widetilde E$ is $\mu_Y$-exceptional and $\operatorname{mult}_{\widetilde E}F>0$. Since $N(\mu_Y^*K_Y)$ is effective, we get $\nu_E(K_X)\ge\operatorname{mult}_{\widetilde E}F>0$, which is (i). If $E$ is not contracted, then $\widetilde E$ is also the strict transform of $E'$ under $\mu_Y$, so by \cref{lem:zariski}.\ref{zar:strict} its coefficient in $N(\mu_Y^*K_Y)$ is $\nu_{E'}(K_Y)$. Since $\widetilde E\not\subset\operatorname{Supp}F$, its coefficient in $F$ is zero, and therefore $\nu_E(K_X)=\nu_{E'}(K_Y)$, which is (ii).
	\end{proof}
	
	We now show that an MMP with scaling contracts all components of the negative part after finitely many steps. For projective varieties, this is \cite[Remark~2.3]{DHP13}.
	
	\begin{lemma}\label{lem:prefix}
		Let $X$ be a strongly $\Q$-factorial compact K\"ahler variety with terminal singularities and $K_X$ pseudo-effective. Let $\omega$ be a K\"ahler class such that $K_X+\omega$ is K\"ahler, and run a $K_X$-MMP with scaling of $\omega$ $$X \dashrightarrow X_1 \dashrightarrow X_2 \dashrightarrow \cdots \dashrightarrow X_n \dashrightarrow \cdots,$$
		where
		\[
		\lambda_i : = \inf \{t\ge 0\mid K_{X_i}+t\omega_i \text{ is nef }\}
		\]
		and $\omega_i$ is the push forward of $\omega$ on $X_i$. Either this MMP terminates, or after finitely many steps it reaches a strongly $\Q$-factorial compact K\"ahler variety $S$ with terminal singularities such that $N(K_S)=0$.
	\end{lemma}
	
	\begin{proof}
		Suppose the MMP does not terminate, so that $\lambda_i>0$ and $\lambda_i\to0$. In particular the intervals $[\lambda_i,\lambda_{i-1}]$ cover $(0,\lambda_0]$. Let $D$ be a prime divisor with $\nu_D(K_X)>0$ (that is a component of $N(K_X)$). By \cref{lem:zariski}.\ref{zar:limit},
		\begin{equation}\label{eq:limit}
			\lim_{t\downarrow0}\nu_D(K_X+t\omega)=\nu_D(K_X)>0 .
		\end{equation}
		Suppose $D$ is not contracted by $X\dashrightarrow X_i$, and let $t\in[\lambda_i,\lambda_{i-1}]$. On a common resolution $\mu_X\colon W\to X$, $\mu_{X_i}\colon W\to X_i$, the strict transform $D_W = (\mu_X)_*^{-1}D$ is not $\mu_{X_i}$-exceptional. By \cref{prop:mmp-scaling}.\hyperref[prop:mmp-scaling:3]{(3)}, we have
		\[
		\mu_X^*(K_X+t\omega)=\mu_{X_i}^*(K_{X_i}+t\omega_i)+E_i(t),
		\]
		such that $N\bigl(\mu_X^*(K_X+t\omega)\bigr)=E_i(t)$, where $E_i(t)$ is $\mu_{X_i}$-exceptional. Since $D_W$ is not $\mu_{X_i}$-exceptional, $\mult_{D_W}E_i(t)=0$. Hence
		\[
		\nu_D(K_X+t\omega)=\nu_{D_W}\bigl(\mu_X^*(K_X+t\omega)\bigr)=\mult_{D_W}E_i(t)=0 .
		\]
		If $D$ were never contracted, this would hold for all $t\in(0,\lambda_0]$, contradicting \eqref{eq:limit}. As $N(K_X)$ has finitely many components, all of them are contracted after finitely many steps. Let $S$ be the model that we obtained, we have $N(K_S)=0$. Terminal singularities and strong $\Q$-factoriality are preserved by the steps of the MMP.
	\end{proof}

	\begin{corollary}\label{lem:flips}
		Let $X$ be a strongly $\Q$-factorial compact K\"ahler variety with terminal singularities and $K_X$ pseudo-effective and $N(K_X) = 0$. Then every step of a $K_X$-MMP with scaling is a flip.
	\end{corollary}
	
	\begin{proof}
		If some step contracted a divisor, then by \cref{lem:contractedN}.\hyperref[lem:contractedN:contracted]{(i)}, $\nu_E(K_X)>0$, a contradiction. Since $K_X$ is pseudo-effective, the MMP does not end with a Mori fiber space, hence every step of the MMP is flip.
	\end{proof}

	\begin{proposition}\label{lem:smallmodel}
	Under the good divisorial Zariski decomposition \cref{prop:goodzariski}, let $X$ be an $n$-dimensional strongly $\Q$-factorial compact K\"ahler variety with klt singularities, with $K_X$ pseudo-effective and $N(K_X)=0$. Then there exists a strongly $\Q$-factorial compact K\"ahler variety $M$ with klt singularities such that $K_M$ is semiample. Moreover, the induced bimeromorphic map $\varphi\colon X\dashrightarrow M$ is isomorphic in codimension 1.
	\end{proposition}

\begin{proof}
	We divide the proof into several steps. We first take a projective resolution $Y\to X$ and use \cref{prop:goodzariski} to construct a klt pair $(Y,\Delta)$ whose adjoint class $K_Y+\Delta$ has a good divisorial Zariski decomposition. We then run a MMP with scaling, which contracts precisely the divisors exceptional over $X$. The resulting model $M$ has semiample canonical divisor and the induced bimeromorphic map $X\dashrightarrow M$ is an isomorphism in codimension one. Below are the details of the proof.

	\hypertarget{lem:smallmodel:step1}{Step 1}. Let $\mu\colon W\to X$ be a projective log resolution and put $E=\Exc(\mu)$. Since $X$ is klt, we can write $K_W+E=\mu^*K_X+F$, where $F\ge0$ is $\mu$-exceptional. As $K_X$ is pseudo-effective, so is $K_W+E$. \cref{prop:goodzariski} gives a modification $\nu\colon Y\to W$ from a compact K\"ahler manifold $Y$ such that $\Exc(\mu_X)$ is an snc divisor, where $\mu_X:=\mu\circ\nu\colon Y\to X$, and
	\[
	\nu^*(K_W+E)\sim_\Q P+R,\qquad P\ \text{semiample},\qquad R=N\bigl(\nu^*(K_W+E)\bigr).
	\]
	Put $G=\nu^*F\ge0$. Then
	\[
		\mu_X^*K_X+G\sim_\Q P+R,\qquad R=N\bigl(\mu_X^*K_X+G\bigr).
	\]
	The divisor $G$ is $\mu_X$-exceptional, so \cref{lem:zariski}.\ref{zar:exc} gives $R=N(\mu_X^*K_X)+G$. Moreover $N(\mu_X^*K_X)$ is $\mu_X$-exceptional since $N(K_X) = 0$ and apply \cref{lem:zariski}.\ref{zar:strict}. Thus
	\[
		\mu_X^*K_X\sim_\Q P+N_X,\qquad N_X:=N\bigl(\mu_X^*K_X\bigr)=R-G\ge0\ \text{ is $\mu_X$-exceptional}.
	\]

	\hypertarget{lem:smallmodel:step2}{Step 2}. Write $K_Y=\mu_X^*K_X+\sum_k a_kF_k$, where $F_k$ are the $\mu_X$-exceptional divisors and $a_k>-1$. Choose $c_k\in[0,1)$ with $a_k+c_k>0$ and put $\Delta=\sum_k c_kF_k$. Then $(Y,\Delta)$ is klt and $K_Y+\Delta=\mu_X^*K_X+\Gamma$, where $\Gamma=\sum_k(a_k+c_k)F_k$ is $\mu_X$-exceptional with $\Supp\Gamma=\Exc(\mu_X)$. By the decomposition of $\mu_X^*K_X$ in \hyperlink{lem:smallmodel:step1}{Step 1} and \cref{lem:zariski}.\ref{zar:exc},
	\begin{equation}\label{eq:theta}
		K_Y+\Delta\sim_\Q P+\Theta,\qquad \Theta:=N_X+\Gamma=N\bigl(K_Y+\Delta\bigr),\qquad \Supp\Theta=\Exc(\mu_X).
	\end{equation}
	In particular, since $\Supp\Theta=\Exc(\mu_X)\supseteq\Supp\Delta$, the components of $\Theta$ are exactly the prime divisors of $Y$ exceptional over $X$.

	\hypertarget{lem:smallmodel:step3}{Step 3}. We next run a $(K_Y+\Delta + bP)$-MMP with scaling on $Y$. 	Choose $m>0$ such that $mP$ is base point free and Cartier. Fix an integer $b>2nm$.
	By \cref{lem:zariski}.\ref{zar:calc} and \eqref{eq:theta},
	\[
	\begin{aligned}
	(1+b)\Theta
	&=N\bigl((1+b)(K_Y+\Delta)\bigr)\\
	&\le N\bigl(K_Y+\Delta+bP\bigr)+N\bigl(b\Theta\bigr)\\
	&\le N\bigl(K_Y+\Delta+bP\bigr)+b\Theta.
	\end{aligned}
	\]
	Since $P$ is nef, \cref{lem:zariski}.\ref{zar:nef} gives
	\[
	N(K_Y+\Delta+bP)=N(K_Y+\Delta+bP)\le N(K_Y+\Delta)=\Theta.
	\]
	Hence $N(K_Y+\Delta+bP)=\Theta$.
	For convenience, set $H:=K_Y+\Delta+bP\sim_\Q(1+b)P+\Theta$.

	The pair $(Y,\Delta)$ is klt, $P$ is a nef $\Q$-Cartier divisor, and $H$ is pseudo-effective. Apply \cref{prop:mmp-scaling} to run the $H$-MMP with scaling of a K\"ahler class $\omega$ on $Y$:
	\[
	X_0=Y\dashrightarrow X_1\dashrightarrow X_2\dashrightarrow\cdots.
	\]
	For each $i$, let $H_i$ and $\omega_i$ be the push-forwards of $H$ and $\omega$, respectively onto $X_i$, let $\lambda_i:=\inf\{t\ge0\mid H_i+t\omega_i\text{ is nef}\}$ be the scaling threshold, and let $f_i\colon X_i\to Z_i$ be the contraction of the extremal ray $R_i$ defining the $i$-th step. By \cref{prop:mmp-scaling}.\hyperref[prop:mmp-scaling:2]{(2)}, $(H_i+\lambda_i\omega_i)\cdot R_i=0$ and $H_i\cdot R_i<0$.

	\hypertarget{lem:smallmodel:step4}{Step 4}. We claim that each step of this MMP is $P$-trivial. Indeed, if $P_i\cdot R_i>0$ then $(K_{X_i}+\Delta_i)\cdot R_i<0$, so by the cone \cref{thm:kahler-cone} the ray $R_i$ is spanned by a rational curve $C$ with $-(K_{X_i}+\Delta_i)\cdot C\le2n$. Since $mP_i$ is Cartier, $H_i\cdot C\ge b/m-2n>0$, a contradiction.

	Since $mP$ is base point free, let $\phi\colon Y\to V$ be the contraction induced by $|mP|$. Then $V$ is a normal projective variety and $P\sim_\Q\phi^*A$ for an ample $\Q$-Cartier divisor $A$ on $V$ with $mA$ Cartier.

	We next show by induction on $i$ that each $X_i$ admits a morphism $\phi_i\colon X_i\to V$ such that, with $P_i:=\phi_i^*A$ and $\Theta_i$ the push-forward of $\Theta$,
	\begin{equation}\label{eq:Ti}
		H_i\sim_\Q(1+b)P_i+\Theta_i.
	\end{equation}
	For $i=0$, we take $\phi_0=\phi$; then \eqref{eq:Ti} follows from \eqref{eq:theta}. Assume the statement for $i$. Let $C$ be a curve contracted by $f_i$. Its class lies in $R_i$, so $\phi_i^*A\cdot C=P_i\cdot C=0$. Since $A$ is ample, $\phi_i$ contracts $C$. Thus $\phi_i$ is constant on the fibers of $f_i$, and by the rigidity lemma (\cite[Lemma 4.1.13]{BS95}) $\phi_i=\psi_i\circ f_i$ for a morphism $\psi_i\colon Z_i\to V$. Let $f_i^+\colon X_{i+1}\to Z_i$ be the flipped contraction, or the identity if the step is divisorial, and set $\phi_{i+1}=\psi_i\circ f_i^+$. The commutative diagram can be shown as follows:
	\[
	\begin{tikzcd}[column sep=small]
	X_i \arrow[rr, dashed] \arrow[dr, "f_i"'] \arrow[ddr, bend right=30, "\phi_i"'] & & X_{i+1} \arrow[dl, "f_i^+"] \arrow[ddl, bend left=30, "\phi_{i+1}"] \\
	& Z_i \arrow[d, "\psi_i"] & \\
	& V &
	\end{tikzcd}
	\]
	Then $P_{i+1}=f_i^{+*}\psi_i^*A$ is the push-forward of $P_i=f_i^*\psi_i^*A$, and \eqref{eq:Ti} holds for $i+1$-th step.

	\hypertarget{lem:smallmodel:step5}{Step 5}. We next show that every component of $\Theta$ will be contracted by the MMP. Since $P_i\cdot R_i=0$, \eqref{eq:Ti} gives
	\[
	\Theta_i\cdot R_i=H_i\cdot R_i<0.
	\]
	As $\Theta_i$ is effective, every curve with class in $R_i$ is contained in $\Supp\Theta_i$. In particular, if the step is divisorial, the contracted divisor is a component of $\Theta_i$. The components of $\Theta_i$ are strict transforms of components of $\Theta$, hence exceptional over $X$, so no strict transform of a prime divisor of $X$ is ever contracted.

		We next show that every component $Q$ of $\Theta$ is contracted at some step. We shall prove it by contradiction, suppose $Q$ survives at each step of the MMP. By \hyperlink{lem:smallmodel:step3}{Step 3} of our proof, $N(H)=\Theta$, so $\mult_Q N(H)=\mult_Q\Theta>0$. Applying \cref{lem:zariski}.\ref{zar:neflimit} with $\alpha=\{H\}$ and $\gamma=\omega$ gives
		\[
		\lim_{t\downarrow0}\mult_Q N(H+t\omega)=\mult_Q N(H)>0.
		\]
		Therefore there exists $t_0>0$ such that $\mult_Q N(H+t\omega)>0$ for every $0<t<t_0$. For a common resolution $\mu_{X_0}\colon W\to X_0=Y$, $\mu_{X_i}\colon W\to X_i$ as in \cref{prop:mmp-scaling}.\hyperref[prop:mmp-scaling:3]{(3)}. Let $Q'$ be the strict transform of $Q$ (which exists since we assume that $Q$ is not contracted by any step of the MMP). By \cref{lem:zariski}.\ref{zar:strict},
	\[
	\mult_{Q'}N\bigl(\mu_{X_0}^*(H+t\omega)\bigr)=\mult_Q N\bigl(H+t\omega\bigr).
	\] 
	
	If the MMP terminates at $X_N$, then $\lambda_N=0$ and \cref{prop:mmp-scaling}.\hyperref[prop:mmp-scaling:3]{(3)} at $t=0$ gives $\mult_{Q'}E_N(0)=\mult_Q\Theta>0$. This is impossible since $E_N(0)$ is $\mu_{X_N}$-exceptional. On the other hand, $Q'$ maps onto the strict transform of $Q$, which is a prime divisor of $X_N$. Therefore, $\mult_{Q'}E_N(0)=0$. On the other hand, if the MMP does not terminate, then all $\lambda_i$ are positive and $\lambda_i\to0$ by \cref{prop:mmp-scaling}.\hyperref[prop:mmp-scaling:4]{(4)}. Choosing $i$ with $\lambda_{i-1}<t_0$, \cref{prop:mmp-scaling}.\hyperref[prop:mmp-scaling:3]{(3)} at $t=\lambda_{i-1}$ gives
	\[
	\mult_{Q'}N\bigl(\mu_{X_0}^*(H+t\omega)\bigr)=\mult_{Q'}E_i(t)=0,
	\]
	a contradiction.

	Therefore, after finitely many steps, every component of $\Theta$ has been contracted. And we denoted the resulting model by $M$. So that $\Theta_M=0$ and $H_M\sim_\Q(1+b)P_M$ is nef.

	\hypertarget{lem:smallmodel:step6}{Step 6}. By \cref{prop:mmp-scaling}.\hyperref[prop:mmp-scaling:2]{(2)}, $M$ is a strongly $\Q$-factorial compact K\"ahler variety. Since the MMP $Y \dashrightarrow M$ contracts exactly the components of $\Theta$ (where $\Supp \Theta =\Exc(\mu_X)$), so the bimeromorphic map $\varphi\colon X\dashrightarrow M$ neither contracts nor extracts a divisor, and is therefore an isomorphism in codimension one
	\begin{center}
		\begin{tikzcd}
	Y && M & \\
	& X && W.
	\arrow[dashed, from=1-1, to=1-3]
	\arrow["{\mu_X}"', from=1-1, to=2-2]
	\arrow["\varphi"', dashed, from=2-2, to=1-3]
	\arrow["{\mu_M}"', from=2-4, to=1-3]
	\arrow["{\mu_X}", from=2-4, to=2-2]
\end{tikzcd}
	\end{center}
	
	Let $P=\phi^*A$, $P_M=\phi_M^*A$. The pair $(M,b\bar P_M)$ is gklt with nef part $bP_M$ descending to $M$, so $M$ is klt. Finally, since $\Supp\Delta\subseteq\Supp\Theta$, we have $\Supp \Delta_M\subset\Supp \Theta_M$ hence $\Delta_M = 0$. So that $K_M+bP_M=H_M\sim_\Q(1+b)P_M$, and then
	\[
	K_M\sim_\Q P_M=\phi_M^*A
	\]
	is semiample.
\end{proof}

	We next show that, after a small modification, $K_X+s\omega$ stays nef for all sufficiently small $s\ge0$. The proof uses \cref{lem:smallmodel,prop:relmodel} and the cone \cref{thm:kahler-cone}.
	
\begin{proposition}\label{lem:interval}
Under the good divisorial Zariski decomposition \cref{prop:goodzariski}, let $X$ be a strongly $\Q$-factorial compact K\"ahler variety with terminal singularities such that $K_X$ is pseudo-effective and $N(K_X)=0$, and let $\omega$ be a K\"ahler class on $X$. Then there exist a strongly $\Q$-factorial compact K\"ahler variety $Y$, a bimeromorphic map $f\colon X\dashrightarrow Y$ which is an isomorphism in codimension one, and a real number $\delta>0$ such that
	\begin{equation}\label{eq:interval}
		K_Y+s\omega_Y\ \text{is nef for all}\ 0\le s\le\delta,
	\end{equation}
	where $\omega_Y:=f_*\omega$. 
\end{proposition}
	
\begin{proof}
	The strategy is as follows. By \cref{lem:smallmodel} we may replace $X$ by a small model whose canonical class is semiample, hence it is the pullback of an ample divisor on a projective variety $Z$. We then run an MMP over $Z$ for a certain adjoint class. This class has vanishing negative part of divisorial Zariski decomposition, so the MMP contracts no divisor and the output $Y$ is still isomorphic to $X$ in codimension one. Finally, by cone theorem, $K_Y+s\omega_Y$ is nef for $s$ in some interval.

	\hypertarget{lem:interval:step1}{Step 1}. Since $X$ is terminal, hence klt. \cref{lem:smallmodel} gives a small bimeromorphic map $\psi:X \dashrightarrow M$ where $K_M$ is semiample. Choose $m>0$ with $mK_M$ Cartier and globally generated. Let $g\colon M\to Z$ be the induced contraction, such that $Z$ is normal and projective, and
	\begin{equation}\label{eq:KMpullback}
		mK_M\sim g^*H,\qquad H\ \text{an ample Cartier divisor on } Z .
	\end{equation}
	Choose a common resolution $\mu_X\colon W\to X$, $\mu_M\colon W\to M$ of $\psi$ and put $\omega_M=(\mu_M)_*\mu_X^*\omega$. The class $\mu_M^*\omega_M-\mu_X^*\omega$ pushes forward to zero under $\mu_M$. By Demailly's support theorem (\cref{thm:2ndsupp}), it is represented by an $\R$-divisor $D$ exceptional over $M$. Since $\psi$ is an isomorphism in codimension one, $D$ is also exceptional over $X$. For a curve $C$ contracted by $\mu_M$,
	\[
	(-D)\cdot C=\mu_X^*\omega\cdot C\ge0,
	\] so $D\ge0$ by negativity \cref{lem:negativity}. Therefore
	\[
	\mu_M^*\omega_M=\mu_X^*\omega+D
	\]
	is big, so $\omega_M$ is modified big. By \cref{lem:zariski} , then the negative part satisfies
	\begin{equation}\label{equ:neg}
		N(\mu_M^*\omega_M)=N(\mu_X^*\omega)+D=D,
	\end{equation}
	since $\mu_X^*\omega$ is nef and $D$ is effective exceptional over $M$.

	\hypertarget{lem:interval:step2}{Step 2}. Run a relative MMP over $Z$. For $0<t\ll 1$, consider the generalized pair $(M,t\overline{\mu_X^*\omega})$; its trace on $M$ is $t(\mu_M)_*\mu_X^*\omega=t\omega_M$. Set $\alpha:=K_M+t\omega_M$. Since $M$ is klt, the pair $(M, t \overline{\mu_X^* \omega})$ is generalized klt. Then, by \Cref{prop:relmodel}, there exists a bimeromorphic contraction $\varphi\colon M\dashrightarrow Y$ over $Z$, with $g_Y\colon Y\to Z$ and $\omega_Y=\varphi_*\omega_M$. Such that $(Y,t\omega_Y)$ is generalized klt, and
	$\alpha_Y:= K_Y+t\omega_Y$ is nef over $Z$
	\begin{center}
		\begin{tikzcd}
			X & M & Y \\
			& Z
			\arrow["\psi", dashed, from=1-1, to=1-2]
			\arrow["\varphi", dashed, from=1-2, to=1-3]
			\arrow["g"', from=1-2, to=2-2]
			\arrow["{g_Y}", from=1-3, to=2-2]
		\end{tikzcd}.
	\end{center}
	We next show that $\varphi$ is small. The class $K_M$ is nef, and $N(\omega_M)=(\mu_M)_*N(\mu_M^*\omega_M)=(\mu_M)_*D=0$ by \eqref{equ:neg}. Hence for $\alpha=K_M+t\omega_M$, by \cref{lem:zariski} $N(\alpha)=0$. By \cite[Theorem~A.11]{das2026mmpgeneralizedpairskahler}, every divisor contracted by $\varphi$ lies in $\Supp N(\alpha)=\emptyset$. Thus $\varphi$ is small, so $Y$ is isomorphic to $M$, and hence to $X$, in codimension one.

	\hypertarget{lem:interval:step3}{Step 3}. We next show that $K_Y+s\omega_Y$ is nef for all $s$ on some interval. The map $\varphi$ is defined over $Z$ and isomorphism in codimension 1, hence pushing \eqref{eq:KMpullback} forward gives $mK_Y\sim g_Y^*H$. Hence $K_Y$ is nef, and for every curve $C\subset Y$,
	\begin{equation}\label{eq:dichotomy}
		K_Y\cdot C=0\ \text{if $g_Y(C)$ is a point},\qquad K_Y\cdot C\ge\tfrac1m\ \text{otherwise}.
	\end{equation}
	To see this, since $H$ is ample Cartier divisor on $Z$, so that if $g_{Y,*}(C) \ne 0$, then
	\[
	mK_Y \cdot C= g_{Y,*}(C) \cdot H \ge 1.
	\]
	Fix an integer $b>2dm$. We claim that $\alpha_Y+bK_Y=(1+b)K_Y+t\omega_Y$ is nef. By the cone \cref{thm:kahler-cone},
\[
\overline{\mathrm{NA}}(Y)=\overline{\mathrm{NA}}(Y)_{\alpha_Y\ge0}+\sum_i\mathbb R_{\ge0}[C_i],
\qquad 0<-\alpha_Y\cdot C_i\le 2d,
\]
with each $C_i$ a rational curve. Since $K_Y$ is nef, $\alpha_Y+bK_Y$ is nonnegative on the first summand. Since $\alpha_Y$ is nef over $Z$, no $C_i$ is contracted by $g_Y$, so $K_Y\cdot C_i\ge 1/m$ by \eqref{eq:dichotomy}, and
\[
(\alpha_Y+bK_Y)\cdot C_i\ \ge\ -2d+\frac bm\ >\ 0 .
\]
We next put $\delta=t/(1+b)$. For $0\le s\le\delta$,
	\[
	K_Y+s\omega_Y=\Bigl(1-\frac{s(1+b)}{t}\Bigr)K_Y+\frac st\bigl((1+b)K_Y+t\omega_Y\bigr),
	\]
	is sum of two nef classes, hence $K_Y+s\omega_Y$ is nef. This proves \eqref{eq:interval} for $f=\varphi\circ\psi$.
\end{proof}\label{subsec:termination}

	We can now finish the proof of the main theorem of this section:
	\begin{proof}[Proof of \cref{thm:mmp}]
		By \cref{lem:prefix}, either the $K_X$-MMP with scaling terminates, or it reaches a K\"ahler variety $S$ with $K_S$ pseudo-effective and  $N(K_S)=0$. So we may assume the latter.
		
		Choose a K\"ahler class $\omega$ on $S$ such that $K_S+\omega$ is K\"ahler. By \cref{prop:mmp-scaling}, there is a $K_S$-MMP with scaling of $\omega$,
		\[
		S=S_0\dashrightarrow S_1 \dashrightarrow S_2\dashrightarrow\cdots,
		\]
		with thresholds $\lambda_i = \inf \{t\ge 0 \mid K_{S_i} +t \omega_i \text{ is nef}\}$ such that
		 \[
		\lambda_0\ge\lambda_1\ge\lambda_2\ge\cdots>0,
		\] and $K_{S_i}+s\omega_i$ is nef for every $s\in[\lambda_i,\lambda_{i-1}]$, where $\omega_i$ is the push-forward of $\omega$ from $S$ to $S_i$. If this MMP is infinite, then $\lambda_i\to0$. By \cref{lem:flips}, every step is a flip, so every $S_i$ is isomorphic to $S$ in codimension one. By \cref{lem:interval}, there exist a strongly $\Q$-factorial compact K\"ahler variety $Y$, with $S\dashrightarrow Y$ is isomorphism in codimension one, and a constant $\delta>0$ such that $\omega_Y$ is the push-forward of $\omega$ from $S$ to $Y$ and $K_Y+s\omega_Y$ is nef for every $s\in[0,\delta]$
				\[
\begin{tikzcd}
	W & {S_i} \\
	Y & S
	\arrow["{\mu_{S_i}}", from=1-1, to=1-2]
	\arrow["{\mu_Y}"', from=1-1, to=2-1]
	\arrow["{\mu_S}"{description}, from=1-1, to=2-2]
	\arrow[dashed, from=2-2, to=1-2]
	\arrow[dashed, from=2-2, to=2-1].
\end{tikzcd}
		\]
		
		Suppose this MMP does not terminate. Then $\lambda_i>0$ and $\lambda_i\to0$, so the thresholds decrease strictly infinitely often and there is an index $i$ with
		$$0<\lambda_i<\lambda_{i-1}\le\delta .$$
		Choose a smooth common resolution $W$ of $S\dashrightarrow S_i$ and $S\dashrightarrow Y$, with bimeromorphic modification $\mu_S\colon W\to S$, $\mu_{S_i}\colon W\to S_i$, and $\mu_Y\colon W\to Y$.
		Since $S_i$ and $Y$ are isomorphic in codimension one, the prime divisors on $W$ exceptional over $S_i$ are exactly those exceptional over $Y$. Choose compatible canonical divisors. The divisor $\mu_{S_i}^*K_{S_i}-\mu_Y^*K_Y$ pushes forward to zero under $\mu_Y$ and is therefore exceptional over both models. The classes $\mu_{S_i}^*\omega_i-\mu_S^*\omega$ and $\mu_Y^*\omega_Y-\mu_S^*\omega$ push forward to zero under $\mu_{S_i}$ and $\mu_Y$, respectively. By Demailly's support theorem (\cref{thm:2ndsupp}), each is represented by an $\R$-divisor exceptional over both models. Subtracting gives the required representative for $\mu_{S_i}^*\omega_i-\mu_Y^*\omega_Y$.
		
		For $t\in\{\lambda_i,\lambda_{i-1}\}$, the class $K_{S_i}+t\omega_i$ is nef, and $K_Y+t\omega_Y$ is nef by \cref{lem:interval}. Then by \cref{lem:twoparam}, \[
		\mu_{S_i}^*K_{S_i}=\mu_Y^*K_Y,
		\] and since $K_Y$ is nef, $K_{S_i}$ is nef. This contradicts $\lambda_i>0$. Hence the MMP terminates in finite steps, say at $M:= S_n$.
		
		Composing the MMP $X \dashrightarrow S$ and $S \dashrightarrow M$ gives a $K_X$-MMP $X\dashrightarrow M$ with $K_M$ nef, and $M$ is strongly $\Q$-factorial compact K\"ahler variety with terminal singularities. By \cref{prop:abundance}, $K_M$ is semiample. On a common resolution $\mu_X\colon W\to X$, $\mu_M\colon W\to M$, we have \[
		\mu_X^*K_X\sim_\Q \mu_M^*K_M+F
		\] with $F\ge0$ exceptional over $M$, and $F$ is positive along every divisor contracted by $X \dashrightarrow M$. So that $X\dashrightarrow M$ is a good minimal model.
	\end{proof}
	
	\section{Invariance of plurigenera}\label{sec:proof}
	
	In this section, we first recall the invariance of plurigenera established by Levine--Koll\'ar--Mori \cite{Levine83,Levine85,KM92} when the central fiber has a semiample canonical divisor. We then combine this result with \cref{thm:mmp,thm:deformation} to prove the main \cref{thm:main}.
	
\subsection{Invariance of plurigenera with semiample canonical divisor on the central fiber}
The following invariance of plurigenera result, due to Levine--Koll\'ar--Mori \cite{Levine83,Levine85,KM92}, is a key ingredient in our argument. We include a proof for completeness and for reader's convienence.

\begin{theorem}[{Levine \cite{Levine83,Levine85}, Koll\'ar--Mori \cite[Theorem~12.5.5]{KM92}}]\label{thm:lkm}
Let $f\colon\cX\to\Delta$ be proper and flat. Suppose that $X_0$ has terminal singularities, belongs to Fujiki's class $\mathcal{C}$, and that $\omega_{X_0}^{[m]}$ is invertible and globally generated for some $m>0$. Then there exists an analytic neighborhood $U_m$ of $0$ such that $P_k(X_t)=P_k(X_0)$ for all $k= 1,\ldots, m$ and $t\in U_m$.
\end{theorem}

In \cite{Levine83}, Levine considers a smooth family whose central fibre $X_0$ is a compact complex manifold in the Fujiki's class $\mathcal{C}$ and shows that a smooth divisor $D\in|mK_{X_0}|$ lifts to every infinitesimal deformation of $X_0$. So $P_m$ is locally constant as soon as the general member of $|mK_{X_0}|$ is smooth, for instance when $|mK_{X_0}|$ is base point free. In fact, we pass to the $m$-fold cyclic cover $\pi\colon Y_0\to X_0$ branched along $D$, which is again a compact K\"ahler manifold. Since
\[
	\pi_*\omega_{Y_0}\;\cong\;\bigoplus_{j=1}^{m}\omega_{X_0}^{\otimes j},
\]
the space $H^0(X_0,\omega_{X_0}^{\otimes m})$ is a direct summand of $H^0(Y_0,\omega_{Y_0})$. The Hodge-to-de Rham spectral sequence of $Y_0$ degenerates at $E_1$. Hence holomorphic top forms on $Y_0$ lift to every infinitesimal deformation of $Y_0$. The section defining $D$ therefore lifts as well. In \cite{Levine85}, the smoothness assumption on $D$ is weakened to allow mild singularities.

Koll\'ar--Mori \cite[Theorem~12.5.5]{KM92} extend Levine's argument to singular fibres. They assume that $X_0$ is projective with log terminal singularities and smooth in codimension two, such that $\omega_{X_0}^{[r]}$ is invertible and globally generated. Terminal singularities satisfy the condition. A general section of $\omega_{X_0}^{[r]}$ extends and defines a cyclic cover $\mathcal{Y}\to \mathcal{X}$. In place of the $E_1$-degeneration used by Levine, Koll\'ar--Mori with the DuBois property. They show that
\[
H^i(Y_0,\mathbb{C})\longrightarrow H^i(Y_0,\cO_{Y_0}),
\]is surjective.

Koll\'ar \cite[comments after Conjecture~17]{Kol21} observes that the same argument applies when $X_0$ is in Fujiki's class $\mathcal{C}$. We use this form of the theorem.

Before proving the result, we recall some preliminaries. Let $f\colon\cX\to\Delta$ be the setting in \cref{thm:lkm}. Write $t$ for the coordinate on $\Delta$ and set $n=\dim X_0$. We may shrink $\Delta$ around $0$ when is needed. For $k\ge0$, put
\[
A_k=\C[t]/(t^{k+1}),\qquad X^{(k)}=\cX\times_\Delta\operatorname{Specan}A_k,
\]
We denote the total space by $\cX$, so that $X^{(0)}=X_0$. A superscript $k$ always denotes the order of a thickening, while a subscript $t$ denotes the fiber over $t\in\Delta$. For $m\in\mathbb Z$, let $\omega^{[m]}_{\cX/\Delta}$ be the reflexive hull of $\omega^{\otimes m}_{\cX/\Delta}$, and put $\omega^{[m]}_{X^{(k)}}=\omega^{[m]}_{\cX/\Delta}|_{X^{(k)}}$. For each integer $k\ge0$ let $s_k\in H^0(X^{(k)},\omega^{[m]}_{X^{(k)}})$. The \emph{$m$-fold cyclic cover defined by $s_k$} is
\[
\pi\colon Y^{(k)}=\operatorname{Specan}_{X^{(k)}}\bigoplus_{i=0}^{m-1}\omega^{[-i]}_{X^{(k)}}\longrightarrow X^{(k)};
\]
for $\ell\le k$, $Y^{(k)}\times_{A_k}A_\ell$ is the cyclic cover $Y^{(\ell)}$ defined by $s_\ell=s_k|_{X^{(\ell)}}$. A section $s\in H^0(\cX,\omega^{[m]}_{\cX/\Delta})$ defines a cyclic cover $\pi\colon\cY\to\cX$ by the same construction. For each $k\ge0$, its restriction to $X^{(k)}$ is the cyclic cover $Y^{(k)}$ defined by $s_k=s|_{X^{(k)}}$. We write $Y_0:=Y^{(0)}$ for its central fiber, and $Y_t$ for the fiber over $t\in\Delta$. Let $D\subset Y^{(k)}$ be the ramification divisor of $\pi$. We summarize several properties of cyclic covers that will be used below.

\begin{lemma}\label{lem:cyclic}
	Let $k\ge0$ be a positive integer, let $s_k\in H^0(X^{(k)},\omega^{[m]}_{X^{(k)}})$, and let $\pi\colon Y^{(k)}\to X^{(k)}$ be the cyclic cover defined by $s_k$.
	\begin{enumerate}
		\item\label{lem:cyclic:1} The morphism $\pi$ is finite. There is a tautological section $\sigma_k$ of $(\pi^*\omega_{X^{(k)}})^{**}$ with $\sigma_k^N=\pi^*s_k$, vanishing to order one along $D$, and
		\[
		(\pi^*\omega_{X^{(k)}})^{**}=\omega_{Y^{(k)}/A_k}(-(m-1)D).
		\]
		\item\label{lem:cyclic:2} For the cyclic cover $\pi\colon Y_t\to X_t$ of a fiber for $t\in \Delta$, $\pi_*\omega_{Y_t}=\bigoplus_{j=1}^{m}\omega^{[j]}_{X_t}$, and therefore
		\[
		h^n(Y_t,\cO_{Y_t})=h^0(Y_t,\omega_{Y_t})=\sum_{j=1}^{m}P_j(X_t).
		\]
		\item\label{lem:cyclic:3} If $s_0\in H^0(X_0,\omega^{[m]}_{X_0})$ is general, then $Y_0$ is smooth over $X_0\setminus\operatorname{Sing}X_0$ and has canonical singularities.
		\item There exists a sheaf homomorphism $\iota : \Omega^1_{U/A} \to \pi^* \Omega^1_{W/A_n}((m-1)D)$ which is the inverse of the cotangent map $d\pi^*: \pi^* \Omega^1_{W/A_n} \to \Omega^1_{U/A}$ away from the ramification divisor $D$.
	\end{enumerate}
\end{lemma}

\begin{proof}
	(1) The sheaf $\pi_*\cO_{Y^{(k)}}=\bigoplus_{i=0}^{m-1}\omega^{[-i]}_{X^{(k)}}$ is a coherent $\cO_{X^{(k)}}$-algebra, so $\pi$ is finite \cite[II, 1.7]{Nak04}. The remaining assertions follows directly from \cite[proof of Lemma~12.5.7 and (12.5.8)]{KM92}.

	(2) By duality for the finite morphism $\pi$ \cite[Proposition~5.68]{KM98},
	\[
	\pi_*\omega_{Y_t}=\mathcal Hom(\pi_*\cO_{Y_t},\omega_{X_t})=\bigoplus_{i=0}^{m-1}\mathcal Hom(\omega^{[-i]}_{X_t},\omega_{X_t})=\bigoplus_{i=0}^{m-1}\omega^{[1+i]}_{X_t},
	\]
	where the last equality holds because both sheaves $\mathcal Hom(\omega_{X_t}^{[-i]}, \omega_{X_t})$ and $\omega_{X_t}^{[1+i]}$ are reflexive on the normal variety $X_t$ and agree on its smooth locus. Taking global sections gives $h^0(Y_t,\omega_{Y_t})=\sum_{j=1}^{m}P_j(X_t)$, and $h^0(Y_t,\omega_{Y_t})=h^n(Y_t,\cO_{Y_t})$ by Serre duality.
	
	(3) Since $\omega^{[m]}_{X_0}$ is invertible and globally generated, for general section $s_0$ the $\operatorname{div}(s_0)$ is smooth on $X_0\setminus\operatorname{Sing}X_0$ by Bertini's theorem. Over this locus $Y_0$ is smooth \cite[Proposition~4.1.6]{Laz04}. Indeed, the ramification divisor $D=\operatorname{div}(\sigma_0)$ is a Cartier divisor on $Y_0$ which $\pi$ maps isomorphically onto the smooth divisor $\operatorname{div}(s_0)$. Outside $D$ the section $s_0$ does not vanish, so $\pi$ is \'etale and $Y_0$ is smooth because $X_0$ is. And $Y_0$ admits canonical singularities by \cite[Lemma~12.5.10]{KM92}. 

	(4) We next construct the homomorphism $\iota\colon\Omega^1_{U/A_n} \longrightarrow \pi^*\Omega^1_{W/A_n}((m-1)D)$ whose restriction to $U\setminus D$ is the inverse of $d\pi^*$.

Away from $D$, the morphism $\pi$ is \'etale. Hence,
$d\pi^*: \pi^* \Omega^1_{W/A_n} \to \Omega^1_{U/A_n}$ is an isomorphism. Near a point of $D$, choose
relative coordinates in which the cyclic cover is given by
\[
\pi\colon (y,x_2,\ldots,x_n)\longmapsto (y^m,x_2,\ldots,x_n),
\]
with $D=(y=0)$. In these coordinates,
\[
d\pi^*(dx_1)=m y^{m-1}dy,
\qquad
d\pi^*(dx_a)=dx_a\quad(a\ge2).
\]
The inverse of $d\pi^*$ away from $D$ is therefore
\[
dy\longmapsto\frac{\pi^*dx_1}{m y^{m-1}},
\qquad
dx_a\longmapsto\pi^*dx_a\quad(a\ge2).
\]
Since $\mathcal O_U((m-1)D)$ is locally generated by
$y^{-(m-1)}$ (that is every local section can be written as $fy^{-(m-1)}$ for some local holomorphic function $f$).

Since locally $\Omega_{W / A_n}^1=\bigoplus_{k=1}^n \mathcal{O}_W d x_k $, hence $\pi^* \Omega_{W / A_n}^1=\bigoplus_{k=1}^n \mathcal{O}_U \pi^* d x_k$ and therefore by our previous discussion $$\pi^* \Omega_{W / A_n}^1((m-1) D)=\bigoplus_{k=1}^n \mathcal{O}_U\left(\pi^* d x_k \otimes y^{-(m-1)}\right),\quad \Omega^1_{U/A_n} = \bigoplus_{k=1}^n \mathcal{O}_W dx_k.$$Hence one can define a homomorphism $$\iota: \pi^* \Omega_{W / A_n}^1((m-1) D) \to \Omega^1_{U/A_n},$$which is defined under the local frame as $\operatorname{diag}\left(\frac{1}{m}, y^{m-1}, \ldots, y^{m-1}\right)$, which is holommorphically component wise.
\end{proof}
	
	\begin{lemma}\label{lem:duboisdef}
		Let $s\in H^0(\cX,\omega^{[N]}_{\cX/\Delta})$ with $s_0=s|_{X_0}$ general, and let $\pi\colon\cY\to\cX$ be the cyclic cover it defines. For each $k\ge0$, denote its restriction to $X^{(k)}$ by $Y^{(k)}$.
		\begin{enumerate}
			\item\label{lem:duboisdef:1} $Y_0$ is a compact complex variety in Fujiki class $\mathcal{C}$ with Du Bois singularities. Moreover, for every $q\ge 0$, the natural map
			\[
			H^q(Y_0,\C)\longrightarrow H^q(Y_0,\cO_{Y_0})
			\]
			is surjective.
			\item\label{lem:duboisdef:2} For every finite integer $k\ge0$ and every $q\ge0$, $H^q(Y^{(k)},\cO_{Y^{(k)}})$ is a free $A_k$-module, compatible with base change to $A_\ell$ for $\ell\le k$, and spanned over $A_k$ by the image of $H^q(Y_0,\C)$.
			\item\label{lem:duboisdef:3} For every $q\ge0$, $h^q(Y_t,\cO_{Y_t})=h^q(Y_0,\cO_{Y_0})$ for every $t$ near $0$.
		\end{enumerate}
	\end{lemma}
	
	\begin{proof}
		\hypertarget{lem:duboisdef:proof:1}{(1)} By \cref{lem:cyclic}.\hyperref[lem:cyclic:3]{(3)}, $Y_0$ has canonical singularities, hence Du Bois. It is a compact complex variety in Fujiki's class $\mathcal C$, since it is finite over $X_0$. The Hodge-to-de Rham spectral sequence degenerates at $E_1$, so the natural map
		\[
		H^q(Y_0,\C)\longrightarrow H^q(Y_0,\mathcal{O}_{Y_0})
		\]
		is surjective.
		
		\hypertarget{lem:duboisdef:proof:2}{(2)} $Y^{(k)}$ and $Y_0$ have the same underlying topological space, and the inclusion of constants $\C_{Y_0}\to\cO_{Y_0}$ factors as
\[
\C_{Y_0}\longrightarrow\cO_{Y^{(k)}}\longrightarrow\cO_{Y_0},
\]
where the second map is the restriction to $Y_0$. Therefore $H^q(Y_0,\C)\to H^q(Y_0,\cO_{Y_0})$ factors through the restriction map
\[
r^q\colon H^q(Y^{(k)},\cO_{Y^{(k)}})\longrightarrow H^q(Y_0,\cO_{Y_0}).
\]
Since $H^q(Y_0,\mathbb{C})\to H^q(Y_0,\mathcal{O}_{Y_0})$ is surjective by (1), so is $r^q$ for every $q$.

Now consider the morphism $f_k\colon Y^{(k)}\to\operatorname{Spec}A_k$, which is proper and flat. The higher direct image $R^q(f_k)_*\cO_{Y^{(k)}}$ is the $A_k$-module $H^q(Y^{(k)},\cO_{Y^{(k)}})$, and the base change map is
\[
H^q(Y^{(k)},\cO_{Y^{(k)}})\otimes_{A_k}\C\longrightarrow H^q(Y_0,\cO_{Y_0})
\]
induced by $r^q$, which is surjective for every $q$. By the cohomology and base change theorem \cite{BS76}, each $H^q(Y^{(k)},\cO_{Y^{(k)}})$ is therefore a free $A_k$-module. And it commutes with base change
\[
H^q(Y^{(k)},\cO_{Y^{(k)}})\otimes_{A_k}A_\ell\cong H^q(Y^{(\ell)},\cO_{Y^{(\ell)}})\qquad\text{for all }\ell\le k.
\]
Finally, let $M\subset H^q(Y^{(k)},\cO_{Y^{(k)}})$ be the $A_k$-submodule generated by the image of the map
\[
H^q(Y_0,\C)\to H^q(Y^{(k)},\cO_{Y^{(k)}}),
\]
induced by the inclusion of constants $\C_{Y_0}\to\cO_{Y^{(k)}}$. Consider the commutative square
\[
\begin{tikzcd}[column sep=large]
	H^q(Y_0,\C) \arrow[r] \arrow[d, two heads] 
	& H^q(Y^{(k)},\cO_{Y^{(k)}}) \arrow[d, two heads] \\
	H^q(Y_0,\cO_{Y_0}) 
	& H^q(Y^{(k)},\cO_{Y^{(k)}})\otimes_{A_k}\C \arrow[l, "\cong"']
\end{tikzcd}
\]
where the left vertical map is surjective by (1), the right vertical map is the map that quotient out the maximal ideal $\mathfrak{m} = (t)$, and the bottom map is the case $\ell=0$ of the base change isomorphism above. Since the left vertical map is surjective and the bottom map is an isomorphism, the composite of the top and right maps is surjective. Hence $M$ surjects onto $H^q(Y^{(k)},\cO_{Y^{(k)}})\otimes_{A_k}\C$, and so $M=H^q(Y^{(k)},\cO_{Y^{(k)}})$ by Nakayama's lemma.

\hypertarget{lem:duboisdef:proof:3}{(3)} Put $g=f\circ\pi\colon \mathcal{Y}\to\Delta$; this morphism is proper, and flat. For every integer $\ell\ge0$, the restriction of $\mathcal{Y}$ to $A_\ell$ is the cyclic cover $Y^{(\ell)}$ defined by $s|_{X^{(\ell)}}$. By \hyperlink{lem:duboisdef:proof:2}{(2)}, the restriction map factors as
\[
H^q(Y^{(\ell+1)},\cO_{Y^{(\ell+1)}})\longrightarrow H^q(Y^{(\ell+1)},\cO_{Y^{(\ell+1)}})\otimes_{A_{\ell+1}}A_\ell\xrightarrow{\ \cong\ } H^q(Y^{(\ell)},\cO_{Y^{(\ell)}}),
\]
where the first map is the natural quotient map. In particular it is surjective for all $\ell$ and $q$, so the inverse limit $\varprojlim_\ell H^q(Y^{(\ell)},\cO_{Y^{(\ell)}})$ surjects onto $H^q(Y_0,\cO_{Y_0})$. By Grauert's comparison theorem \cite[Chapter~VI]{BS76}, this inverse limit is the completion of the stalk $(R^qg_*\cO_{\mathcal{Y}})_0$. Since a finitely generated module and its completion have the same reduction modulo the maximal ideal, the base change map
\[
R^qg_*\cO_{\mathcal{Y}}\otimes\C(0)\longrightarrow H^q(Y_0,\cO_{Y_0})
\]
is surjective for every $q$. By the cohomology and base change theorem \cite[Chapter~III, \S3]{BS76}, every $R^qg_*\cO_{\mathcal{Y}}$ is locally free in a neighbourhood of $0$, and
\[
R^qg_*\cO_{\mathcal{Y}}\otimes\C(t)\cong H^q(Y_t,\cO_{Y_t})
\]
for $t$ in this neighbourhood. Hence $h^q(Y_t,\cO_{Y_t})$ equals the rank of $R^qg_*\cO_{\mathcal{Y}}$, which is $h^q(Y_0,\cO_{Y_0})$.	
	\end{proof}
	
Since $\omega^{[m]}_{\cX/\Delta}$ is invertible and flat over
$\Delta$, for every $k\ge0$ there is an exact sequence
\begingroup
\renewcommand{\theHequation}{thickening}
\begin{equation}
0\longrightarrow\omega^{[m]}_{X_0}
\xrightarrow{\ t^{k+1}\ }
\omega^{[m]}_{X^{(k+1)}}
\longrightarrow\omega^{[m]}_{X^{(k)}}\longrightarrow0.
\tag{1}\label{eq:thickening}
\end{equation}
\endgroup
For a section $s_k\in H^0(X^{(k)},\omega^{[m]}_{X^{(k)}})$, write
$o(s_k)\in H^1(X_0,\omega^{[m]}_{X_0})$ for its image under the
connecting map $o: H^0(X^{(k)},\omega^{[m]}_{X^{(k)}}) \to H^1(X_0,\omega^{[m]}_{X_0})$. Thus $s_k$ lifts to $X^{(k+1)}$ if and only if
$o(s_k)=0$. The technical heart of the proof lies in the following proposition which shows unobstructness of lifting section infinitesimally.

\begin{proposition}\label{prop:lifting}
Let $k\ge0$ be an integer, and let
$s_0\in H^0(X_0,\omega^{[m]}_{X_0})$ be general. If $s_0$ lifts to a
section of $\omega^{[m]}_{X^{(k)}}$, then it lifts to a section of
$\omega^{[m]}_{X^{(k+1)}}$.
\end{proposition}

\begin{proof}
Choose a lift $s_k$ and let
$\pi\colon Y^{(k)}\to X^{(k)}$ be the cyclic cover it defines, with
tautological section $\sigma_k$. Let $W\subset X^{(k)}$ be the open
subspace whose underlying set is
$X_0\setminus\operatorname{Sing}X_0$, put $U=\pi^{-1}(W)$, and set
$S=Y^{(k)}\setminus U$.

Since $X_0$ admits terminal singularities, the singular locus
of $X_0$ has codimension at least three, and $X_0$ is
Cohen--Macaulay. Thus $X^{(k)}$ is Cohen--Macaulay and smooth over
$A_k$ along $W$. By \cref{lem:cyclic}.\hyperref[lem:cyclic:3]{(3)},
the central fiber $Y_0$ has canonical singularities, so it is
Cohen--Macaulay as well. It follows that $Y^{(k)}$ is
Cohen--Macaulay, that $U$ is smooth over $A_k$, and that
$\operatorname{codim}_{Y^{(k)}}S\ge3$.

Step I. We first construct the obstruction class on the cyclic cover. 
Let
\[
\rho\in
\operatorname{Ext}^1_W(\Omega^1_{W/A_k},\cO_W)
\]
be the Kodaira--Spencer class of the extension
$X^{(k+1)}|_W$. Since $\pi|_U\colon U\to W$ is finite flat, pulling
back extensions gives a class
\[
\pi^*\rho\in
\operatorname{Ext}^1_U(\pi^*\Omega^1_{W/A_k},\cO_U).
\]

The tautological section is a section of
$\omega_{U/A_k}(-(m-1)D)$  by \cref{lem:cyclic}. Multiplication by $\sigma_k\in H^0\left(U, \omega_{U / A_k}(-(m-1) D)\right)$, induces a sheaf homomorphism $m_{\sigma_k}:\cO_U \to  \omega_{U/A_k}(-(m-1)D)$ and hence sends $\pi^*\rho$ to
\[
\operatorname{Ext}^1_U\!
\left(\pi^*\Omega^1_{W/A_k},
      \omega_{U/A_k}(-(m-1)D)\right).
\]
After twisting both arguments by $\cO_U((m-1)D)$, regard this as a
class in
\[
\operatorname{Ext}^1_U\!
\left(\pi^*\Omega^1_{W/A_k}((m-1)D),
      \omega_{U/A_k}\right). 
\]
Let
\[
\eta_U\in
\operatorname{Ext}^1_U\!\left(
\pi^*\Omega^1_{W/A_k}((m-1)D),\omega_{U/A_k}\right)
\]
be the class obtained above by pushing forward along multiplication by
$\sigma_k$ and then twisting by $\cO_U((m-1)D)$.

Note that the cyclic cover will induce a natural homomorphism $\iota: \Omega^1_{U/A_k}\to 
\pi^*\Omega^1_{W/A_k}((m-1)D)$ in \cref{lem:cyclic}. Since $\operatorname{Ext}^1$ is contravariant in its first argument,
$\iota$ induces a map
\[
\iota^*\colon
\operatorname{Ext}^1_U\!\left(
\pi^*\Omega^1_{W/A_k}((m-1)D),\omega_{U/A_k}\right)
\longrightarrow
\operatorname{Ext}^1_U\!\left(
\Omega^1_{U/A_k},\omega_{U/A_k}\right).
\]
Set $\alpha_U=\iota^*(\eta_U)$.

By the extension theorem for $\operatorname{Ext}^1$
\cite[12.5.6]{KM92}, restriction induces an isomorphism
\[
\operatorname{Ext}^1_{Y^{(k)}}\!
\left(\Omega^1_{Y^{(k)}/A_k},\omega_{Y^{(k)}/A_k}\right)
\xrightarrow{\ \cong\ }
\operatorname{Ext}^1_U\!
\left(\Omega^1_{U/A_k},\omega_{U/A_k}\right).
\]
Indeed, the complement $S$ has codimension at least three,
$\Omega^1_{Y^{(k)}/A_k}$ is locally free on $U$ and torsion-free in
codimension one, and the dualizing sheaf
$\omega_{Y^{(k)}/A_k}$ is Cohen--Macaulay. Let
\[
\alpha\in
\operatorname{Ext}^1_{Y^{(k)}}\!
\left(\Omega^1_{Y^{(k)}/A_k},\omega_{Y^{(k)}/A_k}\right)
\]
be the unique class restricting to $\alpha_U$.

Step II. Prove the vanishing of the obstruction class.

The relative differential
$d\colon\cO_{Y^{(k)}}\to\Omega^1_{Y^{(k)}/A_k}$
induces a map
\[
d\colon H^{n-1}(Y^{(k)},\cO_{Y^{(k)}})
\longrightarrow
H^{n-1}(Y^{(k)},\Omega^1_{Y^{(k)}/A_k}).
\]
We show that this map is zero. Let $\underline{A_k}$ denote the constant
sheaf on $Y^{(k)}$. Since $d$ is the relative differential over $A_k$,
the composition
\[
\underline{A_k}\longrightarrow\cO_{Y^{(k)}}
\xrightarrow{\ d\ }\Omega^1_{Y^{(k)}/A_k}
\]
is zero. Thus $d$ vanishes on the image of
$H^{n-1}(Y_0,\mathbb{C})$ in
$H^{n-1}(Y^{(k)},\cO_{Y^{(k)}})$.
By \cref{lem:duboisdef}.\hyperref[lem:duboisdef:2]{(2)},
this image spans $H^{n-1}(Y^{(k)},\cO_{Y^{(k)}})$ over $A_k$.
The map $d$ is $A_k$-linear, so it vanishes everywhere.

Serre duality on the compact Cohen--Macaulay space $Y^{(k)}$
identifies the dual map with
\[
d^\vee\colon
\operatorname{Ext}^1_{Y^{(k)}}\!
\left(\Omega^1_{Y^{(k)}/A_k},\omega_{Y^{(k)}/A_k}\right)
\longrightarrow
H^1(Y^{(k)},\omega_{Y^{(k)}/A_k}).
\]
Since $d=0$, its dual map is also zero. In particular,
$d^\vee(\alpha)=0$.

Multiplication by $\sigma_k^{m-1}$ and the trace for the finite map
$\pi$ therefore give the zero class
\begingroup
\renewcommand{\theHequation}{kmob}
\[
\operatorname{ob}(s_k):=
\operatorname{Tr}\!\left(
-\sigma_k^{m-1}d^\vee(\alpha)\right)
=0
\quad\text{in }H^1(X^{(k)},\omega^{[m]}_{X^{(k)}}).
\tag{2}\label{eq:km-ob}
\]
\endgroup

Step III. We end the proof, by showing that obstruction of lifting the sections vanish. We first give the following claim:
\begin{claim}\label{claim:obs-computation}
On $W$,
\[
  \operatorname{ob}(s_k)|_W
  = -(k+1)t^k\,o(s_k)|_W.
\]
\end{claim}

Assuming this claim for the moment, we will prove it later. Since
$\operatorname{codim}(X^{(k)}\setminus W)\ge3$ and
$\omega^{[m]}_{X^{(k)}}$ is invertible, restriction $H^1 (X^{(k)},\omega_{X^{(k)}}^{[m]})\to H^1 (W, \omega_W^{[m]})$
injective. Thus \eqref{eq:km-ob} implies
\[
t^k o(s_k)=0
\quad\text{in }H^1(X^{(k)},\omega^{[m]}_{X^{(k)}}).
\]

If $k=0$, then $o(s_0)=0$, so $s_0$ lifts to $X^{(1)}$.
Suppose $k\ge1$, and abbreviate
$L_j=\omega^{[m]}_{X^{(j)}}$. Multiplication by $t$ gives a
commutative diagram with exact rows:
\[
\begin{tikzcd}[column sep=large]
0 \arrow[r]
& L_0 \arrow[r,"t^k"] \arrow[d,equal]
& L_k \arrow[r] \arrow[d,"t"]
& L_{k-1} \arrow[r] \arrow[d,"t"]
& 0\\
0 \arrow[r]
& L_0 \arrow[r,"t^{k+1}"]
& L_{k+1} \arrow[r]
& L_k \arrow[r]
& 0.
\end{tikzcd}
\]
The induced long exact cohomology sequences contain the
commutative diagram
\[
\begin{tikzcd}[column sep=small]
H^0(X^{(k)},L_k) \arrow[r] \arrow[d,"t"]
& H^0(X^{(k-1)},L_{k-1}) \arrow[r,"\delta"] \arrow[d,"t"]
& H^1(X_0,L_0) \arrow[r,"t^k"] \arrow[d,equal]
& H^1(X^{(k)},L_k) \arrow[d,"t"]\\
H^0(X^{(k+1)},L_{k+1}) \arrow[r]
& H^0(X^{(k)},L_k) \arrow[r,"o"]
& H^1(X_0,L_0) \arrow[r,"t^{k+1}"]
& H^1(X^{(k+1)},L_{k+1}).
\end{tikzcd}
\]
Since $t^k o(s_k)=0$, exactness of the top row gives
$v\in H^0(X^{(k-1)},L_{k-1})$ with $\delta(v)=o(s_k)$.
Commutativity of the middle square then yields
\[
o(tv)=\delta(v)=o(s_k).
\]
Hence $o(s_k-tv)=0$, so the modified section $s_k-tv$ lifts to
$X^{(k+1)}$. Since $tv$ restricts to zero on $X_0$, this section
still restricts to $s_0$, as required.
\end{proof}

\begin{proof}[Proof of \cref{claim:obs-computation}]
Recall that $\pi:Y^{(k)} \to X^{(k)}$ is a $m$-fold cyclic cover. Let $W\subset X^{(k)}$ be the open subset whose under topological space is $X_0 \setminus \text{Sing}(X_0)$ and $U= \pi^{-1}(W)$. Choose a sufficiently fine open cover $\{W_i\}$ of $W$.
Put $U_i=\pi^{-1}(W_i)$. On corresponding neighbourhoods
in $\mathcal X$, choose holomorphic vector fields $\widetilde{\partial}_i$
lifting $\partial_t$. Thus $\widetilde{\partial}_i(t)=1$.
The differences $\theta_{ij}:=\widetilde{\partial}_i-\widetilde{\partial}_j$ gives cocycle $\{\theta_{ij}\}$
represents the Kodaira--Spencer class
\[
\rho\in H^1(W,T_{W/A_k})
\cong
\operatorname{Ext}^1_W
(\Omega^1_{W/A_k},\mathcal O_W).
\]

Step 1. We find a cocycle representing $\alpha_U$. Since $U$ is smooth over $A_k$, the sheaf
$\Omega^1_{U/A_k}$ is locally free. Hence
\[
\operatorname{Ext}^1_U
(\Omega^1_{U/A_k},\omega_{U/A_k})
\cong
H^1\!\left(
U,\mathcal Hom(\Omega^1_{U/A_k},\omega_{U/A_k})
\right)
\cong H^1(U,\Omega^{n-1}_{U/A_k}).
\]
The second isomorphism comes from the contraction isomorphism
\[
\Phi\colon\Omega^{n-1}_{U/A_k}
\xrightarrow{\ \cong\ }
\mathcal Hom(\Omega^1_{U/A_k},\omega_{U/A_k}),
\qquad
\gamma\longmapsto[\beta\mapsto\beta\wedge\gamma].
\]
We regard the tautological section $\sigma_k$ as a relative top form
through the inclusion
\[
\sigma_k\in
H^0(U,\omega_{U/A_k}(-(m-1)D))
\subset H^0(U,\omega_{U/A_k}).
\]

We follow the construction in the proof of \cref{prop:lifting}.
First, we pull back the Kodaira--Spencer class $\rho$ and
multiply by $\sigma_k$. Twisting by $(m-1)D$ then gives a
cocycle $\{\eta_{ij}\}$ representing $\eta_U$, where
\[
\begin{aligned}
\eta_{ij}\colon
\pi^*\Omega^1_{W/A_k}((m-1)D)
&\longrightarrow\omega_{U/A_k},\\
\xi\otimes h
&\longmapsto
\xi(\pi^*\theta_{ij})\,h\sigma_k.
\end{aligned}
\]
Here $h$ is a local section of $\mathcal O_U((m-1)D)$.

Recall the homomorphism constructed in \cref{lem:cyclic}:
\[
\iota\colon\Omega^1_{U/A_k}
\longrightarrow
\pi^*\Omega^1_{W/A_k}((m-1)D).
\]
Away from $D$, it is the inverse of the cotangent map
$d\pi^*\colon\pi^*\Omega^1_{W/A_k}\to\Omega^1_{U/A_k}$.
Since $\alpha_U=\iota^*\eta_U$, it is represented by the
cocycle $\{\alpha_{ij}\}$, where
\[
\alpha_{ij}:=\eta_{ij}\circ\iota
\colon\Omega^1_{U/A_k}\longrightarrow\omega_{U/A_k}.
\]
Thus $\gamma_{ij}:=\Phi^{-1}(\alpha_{ij})$ gives a cocycle
of holomorphic $(n-1)$-forms representing $\alpha_U$.

We now describe $\gamma_{ij}$ away from $D$.
Away from the ramification divisor, $\pi$ is \'etale. Thus the tangent map
\[
d\pi\colon T_{U/A_k}
\xrightarrow{\ \cong\ }\pi^*T_{W/A_k}
\]
is an isomorphism. Define
$\widehat\theta_{ij}:=(d\pi)^{-1}(\pi^*\theta_{ij})$ over there.
For a local one-form $\beta$, put $\xi=\iota(\beta)$.
Then
\[
\begin{aligned}
\alpha_{ij}(\beta)
&=\xi(\pi^*\theta_{ij})\,\sigma_k\\
&=(d\pi^*\xi)(\widehat\theta_{ij})\,\sigma_k\\
&=\beta(\widehat\theta_{ij})\,\sigma_k.
\end{aligned}
\]
Since
\[
\beta\wedge\iota_{\widehat\theta_{ij}}\sigma_k
=\beta(\widehat\theta_{ij})\,\sigma_k,
\]
we obtain
\[
\gamma_{ij}
=\iota_{\widehat\theta_{ij}}\sigma_k
\qquad\text{on }(U_i\cap U_j)\setminus D.
\]
The form $\gamma_{ij}$ is holomorphic on the entire overlap, thus $\gamma_{ij}$ extends the
contraction across $D$.

Step 2. We next compute a cocycle representing
$\operatorname{ob}(s_k)|_W$. By definition,
\[
\operatorname{ob}(s_k)|_W
=
\operatorname{Tr}\!\left(
-\sigma_k^{m-1}d^\vee(\alpha)|_U
\right).
\]
Hence $\operatorname{ob}(s_k)|_W$
is represented by $\left\{
\operatorname{Tr}\!\left(
-\sigma_k^{m-1}d\gamma_{ij}
\right)
\right\}$. Away from $D$, Cartan's formula gives
\[
d\gamma_{ij}
=d(\iota_{\widehat\theta_{ij}}\sigma_k)
=\mathcal L_{\widehat\theta_{ij}}\sigma_k.
\]
Here $d\sigma_k=0$ because $\sigma_k$ is a relative top form.
Applying the Lie derivative along $\hat{\theta}_{ij}$ to $\sigma_k^m=\pi^*s_k$ gives
\[
m\sigma_k^{m-1}
\mathcal L_{\widehat\theta_{ij}}\sigma_k
=
\pi^*(\mathcal L_{\theta_{ij}}s_k).
\]
Taking the trace and using $\operatorname{Tr}\circ\pi^*=m$,
we obtain
\[
\operatorname{Tr}\!\left(
-\sigma_k^{m-1}d\gamma_{ij}
\right)
=-\mathcal L_{\theta_{ij}}s_k
\]
away from the branch divisor.

Both sides are holomorphic sections on $W_i\cap W_j$. The equality therefore holds on the entire overlap.
Consequently,
\[
\operatorname{ob}(s_k)|_W
=
\bigl[\{-\mathcal L_{\theta_{ij}}s_k\}\bigr].
\]

Step 3. We end the proof of the claim by comparing this cocycle $\{-\mathcal{L}_{\theta_{ij}}s_k\}$ with the cocyle representation of $o(s_k)$. Put $L_j=\omega_{X^{(j)}}^{[m]}$. Recall the exact sequence
\[
0\longrightarrow L_0
\xrightarrow{\ t^{k+1}\ }
L_{k+1}
\xrightarrow{\ r\ }
L_k
\longrightarrow0.
\]
After refining the cover, choose local lifts
$\widetilde s_i$ of $s_k|_{W_i}$ to $L_{k+1}$.
These lifts are defined on the corresponding open subsets
of $X^{(k+1)}$.
On $W_{ij}:=W_i\cap W_j$, both lifts restrict to
$s_k|_{W_{ij}}$. Exactness therefore gives a unique section
$c_{ij}\in H^0(W_{ij}\cap X_0,L_0)$ such that
\[
\widetilde s_j-\widetilde s_i=t^{k+1}c_{ij}.
\]

The connecting homomorphism is computed by the diagram
\[
\begin{tikzcd}
	0 & {\check C^0(L_0)} & {\check C^0(L_{k+1})} & {\check C^0(L_k)} & 0 \\
	0 & {\check C^1(L_0)} & {\check C^1(L_{k+1})} & {\check C^1(L_k).} & 0
	\arrow[from=1-1, to=1-2]
	\arrow["{t^{k+1}}", from=1-2, to=1-3]
	\arrow["{\check\delta}", from=1-2, to=2-2]
	\arrow["r", from=1-3, to=1-4]
	\arrow["{\check\delta}", from=1-3, to=2-3]
	\arrow[from=1-4, to=1-5]
	\arrow["{\check\delta}", from=1-4, to=2-4]
	\arrow[from=2-1, to=2-2]
	\arrow["{t^{k+1}}", from=2-2, to=2-3]
	\arrow["r", from=2-3, to=2-4]
	\arrow[from=2-4, to=2-5]
\end{tikzcd}
\]
The cochains are taken on the corresponding covers $\{W_i\}$. Start with $\{s_k|_{W_i}\}$ and choose the
lift $\{\widetilde s_i\}$. Its \v{C}ech differential is
\[
(\check\delta\widetilde s)_{ij}
=\widetilde s_j-\widetilde s_i=t^{k+1}c_{ij}.
\]
Its preimage in $\check C^1(L_0)$ is thus $\{c_{ij}\}$.
By diagram chasing, the connecting homomorphism gives
\[
o(s_k)|_{W\cap X_0}=[\{c_{ij}\}].
\]

Choose local extensions of $\widetilde s_i$ to
$\mathcal X$. Put
\[
u_i:=
\left.\mathcal L_{\widetilde{\partial}_i}\widetilde s_i\right|_{X^{(k)}}.
\]
One can check this does not depend on the choice of the extension. 

The extensions $\tilde{s}_i$ agree on $X^{(k)}$.
Their difference is therefore divisible by $t^{k+1}$.
Write
\[
\widetilde s_j-\widetilde s_i=t^{k+1}c_{ij}
\]
Here we also use $c_{ij}$ for
the resulting ambient section. Its restriction to $X_0$
is the section defined above.

Since $\widetilde{\partial}_j(t)=1$, we have
\[
\begin{aligned}
u_j-u_i
&=
\left.
\left(
\mathcal L_{\widetilde{\partial}_j}(\widetilde s_j-\widetilde s_i)
-\mathcal L_{\widetilde{\partial}_i-\widetilde{\partial}_j}\widetilde s_i
\right)\right|_{X^{(k)}}\\
&=
\left.
\left(
(k+1)t^k c_{ij}
+t^{k+1}\mathcal L_{\widetilde{\partial}_j}c_{ij}
\right)\right|_{X^{(k)}}
-\mathcal L_{\theta_{ij}}s_k\\
&=(k+1)t^k c_{ij}
-\mathcal L_{\theta_{ij}}s_k.
\end{aligned}
\]
The term divisible by $t^{k+1}$ vanishes on $X^{(k)}$. Thus
\[
-\mathcal L_{\theta_{ij}}s_k
-\bigl(-(k+1)t^k c_{ij}\bigr)
=u_j-u_i.
\]
The two cocycles differ by the coboundary
$\check\delta\{u_i\}$. They therefore represent the same
cohomology class. We conclude that
\[
\operatorname{ob}(s_k)|_W
=-(k+1)t^k o(s_k)|_W.
\]
Here $t^k o(s_k)|_W$ denotes the restriction to $W$ of
the image of $o(s_k)$ under the map induced by
$t^k\colon L_0\to L_k$.
\end{proof}

We can now finish the proof of theorem of Levine--Koll\'ar--Mori. 
\begin{proof}[Proof of \cref{thm:lkm}]
We divide the proof into three steps.

\medskip
Step 1: Lift sections to every infinitesimal order.
For $k\ge0$, consider the restriction map
\[
\gamma^{(k)}:H^0(X^{(k)},\omega^{[m]}_{X^{(k)}})
\longrightarrow H^0(X_0,\omega^{[m]}_{X_0}).
\]
We prove by induction on $k$ that $\gamma^{(k)}$ is surjective. This is clear for $k=0$.
Suppose it holds for $k$, then every general section then lifts to
$X^{(k)}$. By \cref{prop:lifting} lifts further to $X^{(k+1)}$. Thus $\gamma^{(k+1)}$ has dense image, hence it's surjective. 

\medskip
Step 2: From infinitesimal lifting to an actual lifting in a neighborhood of $X_0$.
By Step~1 $$
\varprojlim_{k\ge 0} H^0\!\left(X^{(k)},\omega^{[m]}_{\cX/\Delta}|_{X^{(k)}}\right)
\twoheadrightarrow
H^0\!\left(X_0,\omega^{[m]}_{\cX/\Delta}|_{X_0}\right).
$$is surjective. Grauert's comparison
theorem, we can identify this inverse limit with completion of
$(f_*\omega^{[m]}_{\cX/\Delta})_0$ that is$$
\varprojlim_{k\ge 0} H^0\!\left(X^{(k)},\omega^{[m]}_{\cX/\Delta}|_{X^{(k)}}\right)
\cong
\widehat{(f_*\omega^{[m]}_{\cX/\Delta})_0}.$$Since completion preserves reduction
modulo the maximal ideal, the restriction map
\[
\left(f_* \omega^{[m]}_{\cX/\Delta}\right)_0 \longrightarrow \frac{\left(f_* \omega^{[m]}_{\cX/\Delta}\right)_0}{\mathfrak{m}_0\left(f_* \omega^{[m]}_{\cX/\Delta}\right)_0} \longrightarrow H^0\left(X_0,\left.\omega^{[m]}_{\cX/\Delta}\right|_{X_0}\right) .
\]
is surjective. After shrinking $\Delta$, choose
$s\in H^0(\cX,\omega^{[m]}_{\cX/\Delta})$ whose restriction $s_0$ is general.

\medskip
Step 3: Prove $P_m (X_t)$ is locally constant.
Let $\pi\colon\mathcal Y\to\cX$ be the $m$-fold cyclic cover defined by
$s$. By \cref{lem:cyclic}.\hyperref[lem:cyclic:2]{(2)} and
\cref{lem:duboisdef}.\hyperref[lem:duboisdef:3]{(3)}, for $t$ near
$0$,
\[
\sum_{j=1}^{m}P_j(X_t)
=h^n(Y_t,\cO_{Y_t})
=h^n(Y_0,\cO_{Y_0})
=\sum_{j=1}^{m}P_j(X_0).
\]
Each $P_j(X_t)$ is upper semicontinuous. Since their sum is constant, $P_j(X_t)$ is constant near $0$ for $1\le j\le m$. 
\end{proof}
	
	\subsection{Proof of the main theorem}
	We can now prove the main theorem of this paper: 
	\begin{theorem}\label{thm:main}
	Under the good divisorial Zariski decomposition \cref{prop:goodzariski}, let  $\pi\colon\cX\to\Delta$ be a smooth family of compact K\"ahler manifolds. Then for every integer $m\ge1$, $P_m(X_t)$ is independent of $t\in\Delta$.
	\end{theorem}
	
\begin{proof}
Fix $m\ge1$ and $t_0\in\Delta$. If $P_m(X_{t_0})=0$,
upper semicontinuity gives $P_m(X_t)=0$ for $t$ near $t_0$.
We may therefore assume that $P_m(X_{t_0})>0$, so that
$K_{X_{t_0}}$ is pseudo-effective.

By \cref{thm:mmp}, there is a finite $K_{X_{t_0}}$-MMP
$X_{t_0}\dashrightarrow M_{t_0}$, where $M_{t_0}$ is a
compact K\"ahler variety with terminal singularities and
semiample canonical divisor. Choose $r>0$ such that
$\omega_{M_{t_0}}^{[r]}$ is invertible and globally generated.
By \cref{thm:deformation,cor:through}, after restricting to
a sufficiently small disk $\Delta$ centered at $t_0$, this MMP
extends to the family and gives a flat proper morphism
$f\colon\mathcal M\to \Delta$ with central fiber $M_{t_0}$.
After shrinking $\Delta$, its fibers have terminal singularities,
are bimeromorphic to the corresponding fibers $X_t$, and satisfy
\[
P_j(M_t)=P_j(X_t)
\qquad\text{for all }j\ge1\text{ and }t\in \Delta.
\]
In particular, every $M_t$ belongs to Fujiki's class
$\mathcal C$.

We next follow the argument of
\cite[Corollary~1.10]{Levine83}.
Choose $r>0$ such that $\omega_{M_{t_0}}^{[r]}$
is invertible and globally generated.
By \cref{thm:lkm},
$P_r(M_t)=P_r(M_{t_0})$ for $t$ near $t_0$.
Thus \cite[Corollary~3.14]{HLR26} implies, after shrinking the neighborhood further,
that $\omega_{M_t}^{[r]}$ is invertible and globally
generated for every $t$ in a connected neighborhood
$U$ of $t_0$. We may also assume that every $M_t$,
$t\in U$, has terminal singularities. Moreover,
each $M_t$ belongs to Fujiki's class $\mathcal C$,
since it is bimeromorphic to the compact K\"ahler
manifold $X_t$.

Choose $N\in r \mathbb{Z}_+$.
For every $s\in U$, the sheaf $\omega_{M_s}^{[N]}$ is invertible
and globally generated. Applying \cref{thm:lkm} with $M_s$ as
the central fiber, we obtain a neighborhood $U_N(s)\subset U$
of $s$ such that
\[
P_j(M_t)=P_j(M_s),
\qquad
1\le j\le N,\quad t\in U_N(s).
\]
Since two non-empty Zariski open subset has non-empty interesection $U_N(t_0)\cap U_N(s)\ne\emptyset$, we obtain
\begin{equation}\label{eq:pluri}
P_j(M_t)=P_j(M_{t_0}),
\qquad
1\le j\le N,\quad t\in U.
\end{equation}
Since $N$ can be chosen arbitrarily large among the positive
multiples of $r$, \eqref{eq:pluri} holds for every $j\ge1$
and every $t\in U$. By \cref{cor:through}, it follows that
\[
P_j(X_t)=P_j(M_t)=P_j(M_{t_0})=P_j(X_{t_0}),
\qquad
j\ge1,\quad t\in U.
\]
This proves that, for every fixed index $m\ge1$,
the function $t\mapsto P_m(X_t)$ is locally constant at every
$t_0\in U$. Since $\Delta$ is connected, $P_m(X_t)$ is
independent of $t\in\Delta$.
\end{proof}

\textbf{Added in proof.} During our polishing of this paper, Z. Hu--K. Liu \cite{HuLiu26} proved Siu's Conjecture \ref{conj:siu} even in the log version, by using Hacon--Wang's termination result \cite{HW26} (under the assumptions of log Iitaka subadditivity and log abundance conjectures).

	\raggedbottom

\end{document}